\documentclass[11pt,a4paper]{article}%
\usepackage{amsfonts}
\usepackage{amsmath}
\usepackage{amssymb}
\usepackage{graphicx}%
\usepackage[latin1]{inputenc}
\providecommand{\U}[1]{\protect\rule{.1in}{.1in}}
\newtheorem{theorem}{Theorem}

\newtheorem{conjecture}[theorem]{Conjecture}
\newtheorem{corollary}[theorem]{Corollary}

\newtheorem{definition}[theorem]{Definition}
\newtheorem{example}[theorem]{Example}

\newtheorem{lemma}[theorem]{Lemma}

\newtheorem{problem}[theorem]{Problem}
\newtheorem{proposition}[theorem]{Proposition}
\newtheorem{remark}[theorem]{Remark}

\newtheorem{openproblem}[theorem]{Open Problem}

\begin{document}

\title{On Wigdersons' approach to the uncertainty principle}
\author{Nuno Costa Dias\textbf{\thanks{nunocdias1@gmail.com\, (Corresponding author)}} 
\and Franz Luef\textbf{\thanks{franz.luef@ntnu.no}}
\and Jo\~{a}o Nuno Prata\textbf{\thanks{joao.prata@mail.telepac.pt }}}
\maketitle

\begin{abstract}
We revisit the uncertainty principle from the point of view suggested by A. Wigderson and Y. Wigderson. This approach is based on a primary uncertainty principle from which one can derive several inequalities expressing the impossibility of a simultaneous sharp localization in time and frequency. Moreover, it requires no specific properties of the Fourier transform and can therefore be easily applied to all operators satisfying the primary uncertainty principle. A. Wigderson and Y. Wigderson also suggested many generalizations to higher dimensions and stated several conjectures which we address in the present paper. We argue that we have to consider a more general primary uncertainty principle to prove the results suggested by the authors. As a by-product we obtain some new inequalities akin to the Cowling-Price uncertainty principle, a generalization of the Heisenberg uncertainty principle, and derive the entropic uncertainty principle from the primary uncertainty principles. 
\end{abstract}

{\bf Keywords}: Uncertainty principles; Fourier Transform; Metaplectic operators.

{\bf AMS Subject Classifications (2020)}: 81S07; 42A38; 26D10; 46E30

\section{Introduction}

In \cite{Wigderson} the authors revisited the uncertainty principle by showing that several inequalities which involve a function and its Fourier transform are in fact valid, for more general operators than the Fourier transform. They called them \textit{k-Hadamard operators}. An operator $A: L^1 (\mathbb{R}^n)\to L^{\infty} (\mathbb{R}^n)$ is said to be k-Hadamard for some $k>0$, if $\|A\|_{1 \to \infty} \leq 1$, and $\|A^{\ast}Af \|_{\infty} \geq k \|f\|_{\infty}$, for all $f$ such that $f,Af \in L^1 (\mathbb{R}^n)$. In particular, the Fourier transform is $1$-Hadamard. Equipped with this concept they proved several results. The starting point were two fundamental inequalities, which they called the \textit{primary uncertainty principle for k-Hadamard operators},
\begin{equation}
\|f\|_1 \|Af\|_1 \geq k \|f\|_{\infty} \|Af\|_{\infty}~,
\label{eq0A}
\end{equation}
whenever $f,Af \in L^1 (\mathbb{R}^n)$, and the \textit{norm uncertainty principle},
\begin{equation}
\|f\|_1 \|Af\|_1 \geq k^{1-\frac{1}{q}} \|f\|_q \|Af\|_q ~,
\label{eq0B}
\end{equation}
for $f,Af \in L^1 (\mathbb{R}^n)$ and $q \in \left[\right.1, \infty \left.\right)$. As a consequence of these inequalities, they proved the following variation of Heisenberg's uncertainty principle:
\begin{equation}
V(f) V(Af) \geq C_q k^{3-2/q} \|f\|_q^2 \|Af\|_q^2~,
\label{eq1}
\end{equation}
for all $f,Af \in  L^1 (\mathbb{R})$, and $C_q=2^{-\frac{10q-8}{q-1}}$, $q \in \left(\right.1, + \infty \left. \right]$. Here $V(f)$ denotes the variance:
\begin{equation}
V(f) = \int_{\mathbb{R}}x^2 ~|f(x)|^2 dx ~.
\label{eq2}
\end{equation}
If $q=2$ and $Af=\mathcal{F}f= \widehat{f}$ is the Fourier transform (and hence $k=1$), we obtain the familiar Heisenberg uncertainty principle:
\begin{equation}
\left(\int_{\mathbb{R}}x^2 ~|f(x)|^2 dx\right) \cdot \left(\int_{\mathbb{R}}\xi^2 ~|\widehat{f}(\xi)|^2 d \xi \right) \geq C \|f\|_2^2 \|\widehat{f} \|_2^2~,
\label{eq3}
\end{equation}
for some constant $C>0$.

In \cite{Wigderson} the authors left as an open problem the generalization of (\ref{eq1}) to dimension $n>1$. The analogous proof of this higher dimensional inequality would rely on the factorization of a certain integral which is divergent when $n>1$. They suggested that the higher dimensional inequalities could be proved if certain alternative auxilary results were valid. 

In this work, we shall address various conjectures stated in \cite{Wigderson}, and solve the higher dimensional case if $A$ satisfies an additional condition akin to the Hausdorff-Young inequality. Notice that several uncertainty principles for functions belonging to functional spaces other than $L^2$ were derived for various modulation spaces in \cite{Dias1}. Some uncertainty principles were obtained for operators other than the Fourier transform. S. Ghobber \cite{Ghobber} proved variants of the Heisenberg uncertainty principle, the Faris local uncertainty principle and the Donoho-Stark uncertainty principle for certain integral operators, which included the Hankel transform (and also the Dunkl transform \cite{Dunkl1}, \cite{Dunkl2}, \cite{DeJeu2}), the generalized Fourier transform and the $\mathcal{G}$-transform. M.F.E. de Jeu \cite{DeJeu1} also obtained an uncertainty principle for integral operators with a bounded kernel and for which a Plancherel Theorem holds. Hausdorff-Young type inequalities were obtained for certain integral operators by B. Russo \cite{Russo} and by C. Oehring and J.A. Cochran \cite{Oehring}.

We will now briefly state our results and postpone the proofs to later sections.

We shall use the following functionals
\begin{equation}
\mathcal{F}_{p,q}^{(A,n)} (f)=\frac{\|f\|_p \|A f \|_p}{\|f\|_q \|Af\|_q} ~,
\label{eqfunctionals1}
\end{equation}
defined for $p,q \in \left[1, \infty\right]$ and $f$ a measurable complex-valued function on $\mathbb{R}^n$, and belonging to some appropriate domain. When $A = \mathcal{F}$ is the Fourier transform, we write, for simplicity $\mathcal{F}_{p,q}^{(n)}$:
\begin{equation}
\mathcal{F}_{p,q}^{(n)} (f)=\frac{\|f\|_p \|\widehat{f} \|_p}{\|f\|_q \|\widehat{f}\|_q} ~.
\label{eqfunctionals1A}
\end{equation}

The thesis of \cite{Wigderson} is that many uncertainty principles can be derived from the primary and norm uncertainty principles (\ref{eq0A},\ref{eq0B}) and, consequently, it is immaterial whether the operator $A$ is the Fourier transform or some other operator (satisfying (\ref{eq0A},\ref{eq0B})). In particular, if $A$ is a $k$-Hadamard operator, it satisfies inequalities (\ref{eq0A},\ref{eq0B}) and concomitantly all uncertainty principles derived from them.

Here we wish to complete this perspective. Our analysis shows that if we want to derive the higher-dimensional versions of (\ref{eq1}) and (\ref{eq3}), or the entropic uncertainty principle (see Theorem \ref{TheoremEntropicUP} below) in this fashion, then (\ref{eq0A},\ref{eq0B}) will not be sufficient. We will have to consider the following norm uncertainty principles (Theorem \ref{TheoremAdcional1}):
\begin{equation}
\mathcal{F}_{p,q}^{(A,n)} (f) \geq C ~,
\label{IntroFinal1}
\end{equation} 
for all $p \in \left[1,2\right]$ and $q \in \left[p,p^{\prime} \right]$, and where $C>0$ is a constant depending only on $p,q,n$. This is a generalization of the primary and norm uncertainty principles (\ref{eq0A},\ref{eq0B}), to which it reduces for $p=1$.

Arbitrary $k$-Hadamard operators may not satisfy (\ref{IntroFinal1}), so we need to impose more restrictions. The additional conditions are stated in Definition \ref{DefinitionHadamardSpecial} and they ensure that both $A$ and its inverse $A^{-1}$ satisfy a Hausdorff-Young type inequality (Proposition \ref{PropositionHausdorffYoung}). We shall call these operators \textit{special} $k$-\textit{Hadamard operators}. Of course the Fourier transform is a special $k$-Hadamard operator, but so are all metaplectic operators. 

By resorting to the norm uncertainty principles (\ref{IntroFinal1}), we were able to derive various uncertainty principles of the form
\begin{equation}
\|~ |x|^{\theta} f \|_p^{\alpha} ~\| ~| \xi|^{\phi} Af \|_q^{\beta} \geq C \|f\|_r^{\alpha} \|A f\|_r^{\beta}~,
\label{IntroFinal2}
\end{equation} 
for various values of $\alpha, \beta, \theta, \phi, p,q,r >0$ (Theorems \ref{TheoremCowlingPriceInfinity}, \ref{TheoremCowlingPriceHadamard}, \ref{TheoremAdicional3}, \ref{TheoremPreCowlingPrice}), the higher-dimensional version of (\ref{eq1}) and (\ref{eq3}) (Theorem \ref{Theorem2}) and the entropic uncertainty principle (Theorem \ref{TheoremEntropicUP}). We also study the general properties of the sets of $k$-Hadamard operators and special $k$-Hadamard operators (Lemmas \ref{Lemma1}, \ref{Lemma2}, \ref{Lemma3} and Propositions \ref{PropositionInclusion1} and \ref{PropositionInclusion2}).

Shortly before we finished the present paper, Y. Tang published a preprint \cite{Tang} which deals with several related problems, and which complements some of our results for the particular case of the Fourier transform. In particular, the paper presents the following variation of the Heisenberg uncertainty principle for every $f\in {\cal S}({\mathbb{R}^d})\backslash \{0\}$: 
\[
V_p(f) V_p(\widehat{f}) \geq C\, \|f\|_p^p \|\widehat{f}\|_p^p~,
\]
where $1 < p< \frac{2d}{d-1}$, $C=C(d,p)>0$ and $V_p(f)= \int |x|^p ~|f(x)|^p dx$ is the $p$-variance. Note that this uncertainty principle is not an extension of (\ref{eq1}) to the case $q\not=2$.

In the same paper, Y. Tang also obtains the Cowling-Price uncertainty principles in arbitrary dimension using the Wigdersons' approach, and including the endpoints. Note that we have also proved these inequalities in higher-dimension in \cite{Dias1} using a different approach, but our results do not include the endpoints.  

This paper is organised as follows: In Section 2.1 we present the main results of the paper, and in Section 2.2 we discuss some further results, related to the same general topic. In Section 2.3 we investigate the general properties of the sets of $k$-Hadamard and special $k$-Hadamard operators, and in Section 2.4 we state some interesting open problems. The proofs of all results are given in Section 3. The Appendix provides several examples that are useful for the proofs.

\section{Summary of Results}

\subsection{Main results}

Among various open problems, the authors stated the following conjecture (Conjecture 4.13 in \cite{Wigderson}).

\begin{conjecture}\label{Conjecture4.13}
Let $2 \neq p \in \left(\right.1, \infty \left. \right]$, and define the functional $\mathcal{F}_{p,2}^{(1)}: \mathcal{S} (\mathbb{R}) \backslash \left\{0 \right\} \to \mathbb{R}_{>0}$ by
\begin{equation}
\mathcal{F}_{p,2}^{(1)} (f)= \frac{\|f\|_p  \|\widehat{f}\|_p}{\|f\|_2  \|\widehat{f}\|_2}= \frac{\|f\|_p  \|\widehat{f}\|_p}{\|f\|_2^2}~.
\label{eqUP1}
\end{equation}
Then the image of $\mathcal{F}_{p,2}^{(1)}$ is all of $\mathbb{R}_{>0}$.
\end{conjecture}

This conjecture was motivated by the fact that they proved it to be valid, when $p=\infty $. The conjecture was shown in \cite{Huang} and in \cite{Tang} to be valid for $p >2$, but false for $1<p<2$. 

Notice that we can re-express the primary and norm uncertainty principles (\ref{eq0A},\ref{eq0B}) as:
\begin{equation}
\mathcal{F}_{p,q}^{(A,n)} (f) \geq C_{p,q,n}~,
\label{eqNew1}
\end{equation}
for $p=1$, all $q \in \left[1, \infty \right]$ and all $f, Af \in L^1 (\mathbb{R}^n)$. Here $C_{p,q,n}$ is some positive constant, which can only depend on $n$ and $q$. 

It then seems natural to ask whether inequalities (\ref{eqNew1}), which hold for the functional (\ref{eqfunctionals1}) with $p=1$, could also be valid for other values of $p$ and $q$ with some constant $C_{p,q,n}>0$. It is clear (from the case where $A= \mathcal{F}$) that, for some values of the parameters, there are lower positive bounds (e.g. $p=1$, $q \in \left[0, \infty \right]$), while for others the functional can get arbitrarily close to zero (for instance $p>2$, $q=2$). This question has been completely solved in \cite{Tang} for $A= \mathcal{F}$. 

We start by extending the result of \cite{Tang} to a subset of $k$-Hadamard operators (which includes the Fourier transform) and arbitrary dimension $n$.

We shall need the Hausdorff-Young inequality to prove our results. We were unable to prove a Hausdorff-Young inequality for a general $k$-Hadamard operator. Although we have the estimate $\|Af\|_{\infty } \leq \|f\|_1$, there is no counterpart of Plancherel's identity (or some alternative inequality for the $L^2$ norm). The latter would be a crucial element in order to apply the Riesz-Thorin Theorem. We can nevertheless obtain a Hausdorff-Young inequality for $k$-Hadamard operators, if we posit an additional restriction. A natural option would be a condition as stated in Proposition 3.24 of \cite{Wigderson}. However, for our purposes we shall require something more restrictive.

\begin{definition}\label{DefinitionHadamardSpecial}
Let $A$ be a bounded operator $A:L^1 (\mathbb{R}^n) \to L^{\infty} (\mathbb{R}^n)$, with $\|A\|_{1 \to \infty} \leq 1$, and also a bounded and invertible operator $A:L^2 (\mathbb{R}^n) \to L^2 (\mathbb{R}^n)$. Suppose in addition that  and there exists a constant $k>0$, such that
\begin{equation}
A^{\ast}= k A^{-1} ~,
\label{HausdorffYoung1}
\end{equation}
and that $A^{\ast}$ is a also a bounded operator $A^{\ast}:L^1 (\mathbb{R}^n) \to L^{\infty} (\mathbb{R}^n)$. We shall call such an operator a \textbf{special k-Hadamard operator}.
\end{definition}

\begin{remark}\label{RemarkHausdorffYoung0}
Notice that $A^{\ast}$ should be a mapping from $L^{\infty}$ to $L^1$ since $A$ is from $L^1$ to $L^{\infty}$. The adjoint $L^1$ to $L^{\infty}$ does make sense, given that we assume that $A$ is a bounded operator $L^1 \to L^{\infty}$ and, since $L^1$ is a predual of $L^{\infty}$, we can use this to associate to $A$ its adjoint, which has the same norm.
\end{remark}

\begin{remark}\label{RemarkHausdorffYoung}
To have a flavor of the type of operators that satisfy the conditions of the previous definition, let us recall that the Fourier transform is an example of such an operator, as are all metaplectic operators associated with free symplectic matrices (commonly known as Linear Canonical Transforms) \cite{Dias2,Lerner}. In examples \ref{ExampleMultiplicationPhase}, \ref{ExampleDiff} and \ref{ExampleCompStepFunc} we will develop some methods to construct new examples of operators of this type. 

For our purposes it is immaterial whether an operator is a special Hadamard operator. What matters is that it satisfies the primary uncertainty principles (\ref{eqAdicional1}) stated in Theorem \ref{TheoremAdcional1} below. In Propositions \ref{PropositionInclusion1} and \ref{PropositionInclusion2} we will prove that not all Hadamard operators are special Hadamard operators, and that there are operators which satisfy (\ref{eqAdicional1}) but are not special Hadamard operators. The purpose of Definition \ref{DefinitionHadamardSpecial} is to provide a refinement of the $k$-Hadamard operators that satisfy the primary uncertainty principles (\ref{eqAdicional1}).
\end{remark}

Given definition \ref{DefinitionHadamardSpecial}, we can now deduce a Hausdorff-Young inequality for these operators and their inverses.
\begin{proposition}\label{PropositionHausdorffYoung}
Fix $p \in \left[1, 2 \right]$ and let $p^{\prime}$ be its H\"older dual $p^{\prime}= \frac{p}{p-1}$. Let $A$ be a special $k$-Hadamard operator. Then $A$ is a $k$-Hadamard operator.  Moreover, there exists a constant $C >0$, such that
\begin{equation}
\|Af \|_{p^{\prime}} \leq C \|f\|_p ~.
\label{HausdorffYoung2}
\end{equation}
The constant $C$ can only depend on $k,p,n$. 

Conversely, if $Af \in L^p(\mathbb{R}^d)$, then we also have:
\begin{equation}
\|f \|_{p^{\prime}} \leq C \|Af\|_p ~.
\label{HausdorffYoung3}
\end{equation}
\end{proposition}
 
We start by proving the following generalization of the primary and norm uncertainty principles (eqs.(\ref{eq0A},\ref{eq0B})).

\begin{theorem}\label{TheoremAdcional1}
Let $A$ be a special $k$-Hadamard operator, $p \in \left[1,2 \right]$ and $q \in \left[p,p^{\prime}\right]$. Then, there exists a constant $C_{p,q,n}>0$, such that
\begin{equation}
\mathcal{F}_{p,q}^{(A,n)} (f)=\frac{\|f \|_p \|A f\|_p}{ \|f \|_q \|Af\|_q} \geq C_{p,q,n}~,
\label{eqAdicional1}
\end{equation}
for all $f, Af  \in L^p (\mathbb{R}^n)$.
\end{theorem}

\begin{remark}\label{RemarkAdicional1}
If $p=1$ ($p^{\prime}=\infty$), we obtain (\ref{eq0A}) by setting $q= \infty$ and (\ref{eq0B}) with $1\leq q< \infty$. Let us point out that this inequality was derived in \cite{Wigderson} (Theorem 3.25) for matrices. We rederive it here in the infinite dimensional case for completeness.

It is also worth remarking that the extended set of primary uncertainty principles (\ref{eqAdicional1}) are sufficient to prove all other uncertainty principles in this work. Thus they could be our starting point. On the other hand, they appear to be necessary, in the sense that the uncertainty principles (\ref{eq0A}) and (\ref{eq0B}) do not seem to be sufficient to prove the higher dimensional versions of the Heisenberg uncertainty principles (\ref{eq1}) and (\ref{eq3}). They will also permit us to prove the entropic uncertainty principle of Hirschman, Beckner, Bialynicki-Birula and Mycielski \cite{Beckner,Birula,Hirschman}.  
\end{remark}

The inequalities (\ref{eqAdicional1}) can be best understood from the point of view of the embeddings of certain functional spaces. Indeed, let $r \in \left[1, \infty\right]$ and define $\mathcal{B}_A^r (\mathbb{R}^n)= L^r (\mathbb{R}^n) \cap \mathcal{A} L^r (\mathbb{R}^n)$ to be the space of all $f \in \mathcal{S}^{\prime} (\mathbb{R}^n)$ such that:
\begin{equation}
\|f\|_{\mathcal{B}_A^r}:= \|f\|_r + \|A f \|_r < \infty~.
\label{eqEmbeddings1}
\end{equation}  

In the sequel, we shall use the notation $A(f) \gtrsim B(f)$ (resp. $A(f)\lesssim B(f)$), if the there exists a constant $C>0$, such that $A(f) \geq C B(f)$ (resp. $A(f) \leq C B(f)$) for all $f$ in some space. The methods used in this work are not well suited to determine sharp constants, so we will not be concerned with them. 

\begin{theorem}\label{TheoremEmbeddings}
Let $p \in \left[1,2 \right]$, $q \in \left[p, p^{\prime} \right]$ and let $A$ a special $k$-Hadamard operator. Then $\mathcal{B}_A^p (\mathbb{R}^n)$ is continuously embedded in $\mathcal{B}_A^q (\mathbb{R}^n)$,
\begin{equation}
\mathcal{B}_A^p (\mathbb{R}^n)\hookrightarrow \mathcal{B}_A^q (\mathbb{R}^n)~.
\label{eqEmbeddings2}
\end{equation}  
In other words, if $f, A f \in L^p (\mathbb{R}^n)$, then $f, A f \in L^q (\mathbb{R}^n)$, and we have:
\begin{equation}
\|f\|_q + \|Af\|_q \lesssim \|f\|_p + \|Af\|_p~. 
\label{eqEmbeddings3}
\end{equation}  
\end{theorem}

Next we prove two Sobolev type inequalities (Theorems \ref{PropositionProof1} and \ref{TheoremAdcional2}), which will be used to prove various uncertainty principles.

\begin{theorem}\label{PropositionProof1}
Let $u,r,p,\theta,t>0$, be such that
\begin{equation}
p>u,~r \ge u,~1<t \leq \infty, ~\text{and } \theta > n\left(\frac{1}{u}-\frac{1}{p}\right)~.
\label{eqPropositionProof11}
\end{equation}
Then:
\begin{equation}
\|~ |x|^{\theta} f \|_p \gtrsim \|f\|_{\infty} \left(\frac{\|f\|_r}{\|f\|_{\infty}}\right)^{r/u} \left(\frac{\|f\|_r}{\|f\|_{rt}}\right)^{\frac{rt(u\theta p-n(p-u))}{npu(t-1)}}~,
\label{eqPropositionProof12}
\end{equation}
for all $0 \neq f \in L^r (\mathbb{R}^n)\cap L^{\infty} (\mathbb{R}^n)$ if $u\neq r$, and $0 \neq f \in L^r (\mathbb{R}^n)\cap L^{rt} (\mathbb{R}^n)$ if $u=r$.
\end{theorem}

Let us remark that in the previous and in the next theorem, we may have quasi-norms $\| \cdot\|_r$, with $0<r<1$.

\begin{theorem}\label{TheoremAdcional2}
Let $\theta,p>0$ and $q \in \left(\right. p, \infty \left. \right]$. Then
\begin{equation}
\|~|x|^{\theta} f \|_p \gtrsim \left( \frac{\|f\|_p}{\|f\|_q}\right)^{\frac{\theta qp}{n(q-p)}} \|f\|_p ~,
\label{eqTheoremAdcional2.1}
\end{equation}
for all $0 \neq f \in L^p (\mathbb{R}^n)\cap L^q (\mathbb{R}^n)$.
\end{theorem}

The previous inequality includes several familiar inequalities. The fractional Laplacian is defined by:
\begin{equation}
\mathcal{F} \left[\left(- \Delta \right)^s f \right] (\xi)= | \xi|^{2s} \mathcal{F} (f) (\xi)~. 
\label{eqTheoremAdcional2.2}
\end{equation}
If $s=1$, we obtain (minus) the ordinary Laplacian, and if $s= \frac{1}{2}$, we have by Plancherel's Theorem 
\begin{equation}
\|\mathcal{F} \left[\left(- \Delta \right)^{\frac{1}{2}} f \right] \|_2 = \frac{1}{2 \pi}\| \nabla f \|_2~. 
\label{eqTheoremAdcional2.3}
\end{equation}

\begin{corollary}\label{CorollaryAdicional1}
Let $s>0$, $1 \leq p <2$, and set
\begin{equation}
\alpha = \frac{n}{2s} \left(\frac{1}{2}- \frac{1}{p} \right)~.
\label{eqCorollaryAdicional1.1}
\end{equation}
Then we have:
\begin{equation}
\|\left(- \Delta \right)^s f  \|_2 \gtrsim \left(\frac{\|f\|_p}{\|f\|_2 }\right)^{\frac{1}{\alpha}} \|f\|_2 ~,
\label{eqCorollaryAdicional1.2}
\end{equation}
for all $0 \neq f \in L^p (\mathbb{R}^n)\cap L^{2} (\mathbb{R}^n)$.
\end{corollary}

\begin{remark}\label{RemarkAdicional1}
If $s=\frac{1}{2}$, $p=1$, we have $\alpha=- \frac{n}{2}$. Then:
\begin{equation}
\|\nabla f\|_2 \gtrsim  \left(\frac{\|f\|_1}{\|f\|_2 }\right)^{-\frac{2}{n}} \|f\|_2 \Leftrightarrow \|f\|_1^{\frac{2}{n}}  \|\nabla f\|_2 \gtrsim  \|f\|_2^{1+ \frac{2}{n}}~,
\label{eqRemarkAdicional1.1}
\end{equation}
which is Nash's inequality.

\end{remark}

\begin{remark}\label{RemarkAFractional_Laplacian}
Notice that, inspired by (\ref{eqTheoremAdcional2.2}), we may define a fractional $A$-Laplacian for a special $k$-Hadamard operator according to:
\[
\left[(-\Delta_A)^s f \right](x)= A^{-1} \left[|\xi|^{2s}(Af)(\xi) \right]=kA^{\ast} \left[|\xi|^{2s}(Af)(\xi) \right]~.
\]
The previous results then apply directly to the fractional $A$-Laplacian.
\end{remark}

An immediate consequence of Theorems \ref{PropositionProof1} and \ref{TheoremAdcional2} are the following sets of uncertainty principles.

\begin{theorem}\label{TheoremCowlingPriceInfinity}
Let $A$ be a $k$-Hadamard operator, and $\alpha,\beta , \theta,\phi,p,q>0$ be such that:
\begin{equation}
\frac{\theta}{n}+\frac{1}{p}>1~, ~\frac{\phi}{n}+\frac{1}{q}>1~,~ 
\label{eqPreCowlingPriceInfinity1}
\end{equation}
and
\begin{equation}
\alpha \left(\frac{\theta}{n}+\frac{1}{p}\right)=\beta \left(\frac{\phi}{n}+\frac{1}{q} \right)~.
\label{eqPreCowlingPriceInfinity2}
\end{equation}
Then:
\begin{equation}
\|~ |x|^{\theta} f\|_p^{\alpha}~ \|~ |\xi|^{\phi} Af \|_q^{\beta} \gtrsim \|f\|_{\infty}^{\alpha} \|Af\|_{\infty}^{\beta}~,
\label{eqPreCowlingPriceInfinity3}
\end{equation}
for all $f,Af \in L^1 (\mathbb{R}^n)$.
\end{theorem}

\begin{theorem}\label{TheoremCowlingPriceHadamard}
Let $A$ be a $k$-Hadamard operator, and $\alpha,\beta , \theta,\phi,p,q,s$ be such that:
\begin{equation}
s \in \left[1, \infty\right]~, ~ \alpha, \beta >0 ,~p,q >1,~\theta>n \left( 1-\frac{1}{p}\right),~\phi> n\left(1-\frac{1}{q}\right),~ 
\label{eqPreCowlingPriceHadamard1}
\end{equation}
and
\begin{equation}
\alpha \left(\theta-n\left(1-\frac{1}{p}\right)\right)=\beta \left(\phi-n\left(1-\frac{1}{q} \right)\right)~.
\label{eqPreCowlingPriceHadamard2}
\end{equation}
Then:
\begin{equation}
\|~|x|^{\theta}f\|_p^{\alpha}~ \|~|\xi|^{\phi} Af\|_q^{\beta} \gtrsim \|f\|_1^{\alpha} \|Af\|_1^{\beta}~,
\label{eqPreCowlingPriceHadamard3}
\end{equation}
for all $f,Af \in L^1 (\mathbb{R}^n)$. Moreover, if $\alpha=\beta$, then:
\begin{equation}
\|~|x|^{\theta}f\|_p~ \|~|\xi|^{\phi} Af\|_q \gtrsim \|f\|_s \|Af\|_s~.
\label{eqPreCowlingPriceHadamard4}
\end{equation}
\end{theorem}
\begin{remark}\label{RemarkCPH}
Notice that we recover the Heisenberg uncertainty principle (\ref{eq1}) from (\ref{eqPreCowlingPriceHadamard4}) by setting $n=1$, $\theta=\phi=1$ and $p=q=2$, and we also get the additional value $s=1$, which was not considered in \cite{Wigderson}. However, the $n$-dimensional version of (\ref{eq1}) for $n \geq 2$ cannot be obtained from the previous theorem, because the conditions (\ref{eqPreCowlingPriceHadamard1}) for $\theta=1$ and $p=2$ require $n<2$.

Let us also remark that the case $\alpha=\beta$ in Theorem \ref{TheoremCowlingPriceInfinity} coincides with the case $\alpha=\beta$, $s= \infty$ of Theorem \ref{TheoremCowlingPriceHadamard}.
\end{remark}

If we consider only special $k$-Hadamard operators, then we can obtain further inequalities.

\begin{theorem}\label{TheoremAdicional3}
Let $A$ be a special $k$-Hadamard operator, and $\theta>0$, $p \in \left[\right.1,2 \left. \right)$ and $q \in \left[p, p^{\prime}\right]$. Then:
\begin{equation}
\|~|x|^{\theta} f \|_p \|~|\xi|^{\theta} A f \|_p \gtrsim \|f\|_q \|A f\|_q  ,
\label{eqTheoremAdicional3.1}
\end{equation}
for all $f, A f \in L^p (\mathbb{R}^n)$.
\end{theorem}

\begin{theorem}\label{TheoremPreCowlingPrice} 
Let $A$ be a special $k$-Hadamard operator, and let $\alpha,\beta , \theta,\phi,p,q,r$ be such that:
\begin{equation}
1\leq r < 2, ~ \alpha, \beta >0 ,~p,q >r,~\theta> n \left(\frac{1}{r}-\frac{1}{p}\right),~\phi> n \left(\frac{1}{r}-\frac{1}{q}\right),~ 
\label{eqPreCowlingPrice1}
\end{equation}
and
\begin{equation}
\alpha \left[\theta-n \left(\frac{1}{r}-\frac{1}{p}\right) \right]=\beta \left[\phi-n \left(\frac{1}{r}-\frac{1}{q}\right) \right]~.
\label{eqPreCowlingPrice2}
\end{equation}
Then:
\begin{equation}
\|~|x|^{\theta}f\|_p^{\alpha} ~\|~|\xi|^{\phi} A f\|_q^{\beta} \gtrsim \|f\|_r^{\alpha} \|A f\|_r^{\beta}~,
\label{eqPreCowlingPrice3}
\end{equation}
for all $f,A f \in L^r (\mathbb{R}^n)$.
\end{theorem}

\begin{remark}\label{RemarkCowlingPrice1} Let us remark that (\ref{eqPreCowlingPrice3}) is a generalization of (\ref{eqPreCowlingPriceHadamard3}) for $r \neq 1$, but which holds only for special $k$-Hadamard operators.

Notice also that this resembles the Cowling-Price uncertainty principles \cite{Cowling}:
\begin{equation}
\|~|x|^{\theta}f\|_p^{\alpha}~ \|~|\xi|^{\phi} \widehat{f}\|_q^{1-\alpha} \gtrsim \|f\|_2~,
\label{eqCowlingPrice3}
\end{equation}
which hold for all $f\in L^2 (\mathbb{R}^n)$, when  $\alpha,\theta,\phi,p,q>0$ are such that:
\begin{equation}
0<\alpha <1,~p,q \geq 1,~\theta> n\left(\frac{1}{2}-\frac{1}{p}\right),~\phi> n\left(\frac{1}{2}-\frac{1}{q}\right),~ 
\label{eqCowlingPrice1}
\end{equation}
and
\begin{equation}
\alpha \left[\theta-n \left(\frac{1}{2}-\frac{1}{p}\right) \right]=(1-\alpha) \left[\phi-n \left(\frac{1}{2}-\frac{1}{q}\right) \right]~.
\label{eqCowlingPrice2}
\end{equation}
This would be obtained from (\ref{eqPreCowlingPrice3}), by setting $A= \mathcal{F}$, $\beta= 1- \alpha$ and $r=2$. However, the value of $r=2$ is not included in Theorem \ref{TheoremPreCowlingPrice}, so we cannot derive (\ref{eqCowlingPrice3}) in this fashion. On a more positive side, we obtain an inequality which holds not only for the Fourier transform, but also for other special $k$-Hadamard operators and norms $L^r$, $1 \leq r <2$, other than $L^2$. In this respect, the Cowling-Price inequalities (\ref{eqCowlingPrice3}) complete the inequalities (\ref{eqPreCowlingPrice3}) by closing the range of $r$ to $1 \leq r \leq 2$.

\end{remark}

Regarding the higher dimensional generalization of (\ref{eq1}), when $A$ is a special $k$-Hadamard operator, we were able to prove the following theorem (Theorem \ref{Theorem2}). In what follows $V(f)$ is the $n$-dimensional variance:
\begin{equation}
V(f)= \int_{\mathbb{R}^n} |x|^2 ~|f(x)|^2 dx ~,
\label{eq7}
\end{equation}
where this time
\begin{equation}
|x|= \sqrt{x_1^2+x_2^2+ \cdots+x_n^2}~.
\label{eq7A}
\end{equation}
\begin{theorem}\label{Theorem2}
Let $A$ be a special $k$-Hadamard operator. We have:
\begin{equation}
V(f) V(A f) \gtrsim \|f\|_q^2 \|A f\|_q^2~,
\label{eq11}
\end{equation}
for all $f, Af \in L^q(\mathbb{R}^n) \cap L^{q^{\prime}}(\mathbb{R}^n)$, and where:
\begin{equation}
\left\{
\begin{array}{l l}

q \in \left[1, \infty\right] & \text{ , if }n=1\\
\\
q \in \left( 1, \infty \right) & \text{ , if }n=2\\
\\
q \in \left( \frac{2n}{n+2},\frac{2n}{n-2}  \right) & \text{ , if }n \geq 3
\end{array}
\right.
\label{eq11A}
\end{equation}
\end{theorem}

\begin{remark}\label{RemarkHeisenberg}
Notice that Theorem \ref{Theorem2} includes the higher-dimensional Heisenberg uncertainty principle (\ref{eq3}),
\begin{equation}
V(f) V(Af) \gtrsim \|f\|_2^2 \|A f\|_2^2~.
\label{eqHeisenbergHigherDim1}
\end{equation} 
when $q=2$. It is also interesting to remark that, as $n \to \infty$, all that is left is precisely the higher-dimensional Heisenberg uncertainty principle (\ref{eqHeisenbergHigherDim1}).

\end{remark}

This theorem is optimal for the Fourier transform ($A=\mathcal{F}$) in the sense that, when $n \geq 3$, it cannot hold for $q < \frac{2n}{n+2}$ or when $q > \frac{2n}{n-2}$. This is proven in  the next proposition.

\begin{proposition}\label{Proposition3}
Let $n \geq 3$ and $q \in \left[\right.1, \frac{2n}{n+2}\left.\right) \cup \left( \right. \frac{2n}{n-2}, \infty \left.\right]$. Then, there exists a sequence $(f_k)_k$ of functions in $\mathcal{S}(\mathbb{R}^n)$, such that:
\begin{equation}
\frac{V(f_k) V(\widehat{f_k})}{ \|f_k\|_q^2 \|\widehat{f_k}\|_q^2} \to 0~,
\label{eq12}
\end{equation}
as $k \to \infty$.
\end{proposition} 

We next prove the well known entropic uncertainty principle of Hirschman \cite{Beckner,Birula,Hirschman}, by resorting only to the norm uncertainty principles (Theorem \ref{TheoremAdcional1}). As usual we do not obtain optimal constants or optimizers. On the other hand, we generalize it to special $k$-Hadamard operators.
\begin{theorem}\label{TheoremEntropicUP}
Let $A$ be a special $k$-Hadamard operator, such that the constant $C_{p,2,n}$ in (\ref{eqAdicional1}) is differentiable and strictly decreasing with respect to $p$ in some neighborhood $(2-\sigma,2)$ of $2$ ($1 \geq \sigma>0$), and let $f \in \mathcal{S} (\mathbb{R}^n) \backslash \left\{0 \right\}$. Define the Shannon entropies:
\begin{equation}
\left\{
\begin{array}{l}
H\left[f \right]:= - \int_{\mathbb{R}^n}  \left|\frac{f(x)}{\|f\|_2} \right|^2  
\ln \left|\frac{f(x)}{\|f\|_2} \right|^2   dx\\
\\
H\left[A f \right]:= - \int_{\mathbb{R}^n}  \left|\frac{A f(\xi)}{\|A f\|_2} \right|^2  
\ln \left|\frac{Af(\xi)}{\|Af\|_2} \right|^2 d \xi
\end{array}
\right.
\label{eqTheoremEntropicUP1}
\end{equation}
Then, there exists a constant $C_n> -\infty$, such that:
\begin{equation}
H\left[f \right]+ H\left[Af  \right] \geq C_n~,
\label{eqTheoremEntropicUP2}
\end{equation}
for all $f \in \mathcal{S} (\mathbb{R}^n) \backslash \left\{0 \right\}$.
\end{theorem}

Notice that the conditions on the constants $C_{p,2,n}$ are valid for the Fourier transform and for a general metaplectic operator. 

\subsection{Further results}

The following inequalities were also proven in \cite{Wigderson}:
\begin{equation}
\frac{\|f\|_1}{\|f\|_{\infty}} \leq \left(2^5 \frac{V(f)}{\|f\|_{\infty}^2}\right)^{\frac{1}{3}}~,
\label{eq4}
\end{equation}
and
\begin{equation}
\frac{\|f\|_1}{\|f\|_{q}} \leq \left(2^{\frac{5q-4}{q-1}} \frac{V(f)}{\|f\|_{q}^2}\right)^{\frac{q-1}{3q-2}}~,
\label{eq5}
\end{equation}
for $q \in \left(1, \infty \right)$ and $0 \neq f \in L^1 (\mathbb{R})\cap L^{\infty} (\mathbb{R})$. 

The inequality (\ref{eq4}) can be obtained from (\ref{eq5}) by sending $q \to \infty$. By analogy, they proposed the two following open problems in dimension $n>1$.

\begin{openproblem}\label{OpenProblem1}
Let $n \geq 2$. Do there exist constants $\beta_n>0$ and $C_n>0$ (depending only on $n$) such that:
\begin{equation}
\frac{\|f\|_1}{\|f\|_{\infty}} \leq \left(C_n \frac{V(f)}{\|f\|_{\infty}^2}\right)^{\beta_n}~,
\label{eq6}
\end{equation}
for all $f \in L^1(\mathbb{R}^n)\cap L^{\infty}(\mathbb{R}^n)\backslash \left\{0 \right\}$? 
\end{openproblem}

\begin{openproblem}\label{OpenProblem2}
Let $n \geq 2$ and $q \in \left[\right.1, \infty \left. \right)$. Do there exist constants $\beta_{n,q}>0$ and $C_{n,q}>0$ (depending only on $n$ and $q$) such that:
\begin{equation}
\frac{\|f\|_1}{\|f\|_{q}} \leq \left(C_{n,q} \frac{V(f)}{\|f\|_{q}^2}\right)^{\beta_{n,q}}~,
\label{eq8}
\end{equation}
for all $f \in L^1(\mathbb{R}^n)\cap L^{\infty}(\mathbb{R}^n)\backslash \left\{0 \right\}$?
\end{openproblem}

The solution to these problems is negative and it is stated in the following two propostions.

\begin{proposition}\label{Proposition1}
For $n \geq 2$ and all $\beta >0$, there exists a sequence $(f_k)_k$ of functions in $L^1(\mathbb{R}^n)\cap L^{\infty}(\mathbb{R}^n)$, such that
\begin{equation}
\left(\frac{V(f_k)}{\|f_k\|_{\infty}^2}\right)^{\beta} \frac{\|f_k\|_{\infty}}{\|f_k\|_1} \to 0~,
\label{eq9}
\end{equation}
as $k \to \infty$.
\end{proposition}

\begin{proposition}\label{Proposition2}
Let $n \geq 2$, $\beta >0$ and $q  \in \left[\right.1, \infty \left. \right)$ be such that:
\begin{enumerate}
\item $n=2$, or

\item $n>2$ and $q \notin \left(\right. \frac{2n}{n+2},2 \left.\right]$.

\end{enumerate}
Then, there exists a sequence $(f_k)_k$ of functions in $L^1(\mathbb{R}^n)\cap L^{\infty}(\mathbb{R}^n)$, such that
\begin{equation}
\left(\frac{V(f_k)}{\|f_k\|_{q}^2}\right)^{\beta} \frac{\|f_k\|_{q}}{\|f_k\|_1} \to 0~,
\label{eq10}
\end{equation}
as $k \to \infty$.
\end{proposition}

\subsection{Hadamard and special Hadamard operators}

In this section we wish to investigate the sets of Hadamard and special Hadamard operators. 

Given $k>0$, we denote by $\mathcal{H}_k$ the set of $k$-Hadamard operators and by $\mathcal{SH}_k$ the set of special $k$-Hadamard operators. We also define the sets of all Hadamard and special Hadamard operators:
\begin{equation}
\mathcal{H}= \cup_{k > 0} \mathcal{H}_k \hspace{0.5 cm} \text{and} \hspace{0.5 cm} \mathcal{SH}= \cup_{k > 0} \mathcal{SH}_k~.
\label{eqHadSpecHad1}
\end{equation}
Let 
\begin{equation}
p \in \left[1, 2 \right] \hspace{0.5 cm} \text{and} \hspace{0.5 cm} q \in \left[p,p^{\prime}\right]~.
\label{eqHadSpecHad2}
\end{equation}
We define the set $\mathcal{A}_{p,q}$ of all linear operators acting on complex-valued functions defined on $\mathbb{R}^n$, such that there exists a constant $C_{p,q}>0$:
\begin{equation}
\frac{\|Af\|_p \|f\|_p}{\|Af\|_q \|f\|_q} \geq C_{p,q}~, 
\label{eqHadSpecHad3}
\end{equation}
for all $f\neq 0$ such that $f ,Af \in L^p (\mathbb{R})\cap L^q (\mathbb{R})$.

Finally, let
\begin{equation}
\mathcal{A} = \cup_{p \in \left[1, 2 \right]}\cup_{ q \in \left[p,p^{\prime}\right]} \mathcal{A}_{p,q}~.
\label{eqHadSpecHad4}
\end{equation}
Clearly, we have the inclusions:
\begin{eqnarray}
\mathcal{SH}_k \subset \mathcal{H}_k \hspace{1 cm} \label{eqHadSpecHad5A}\\
\mathcal{SH} \subset \mathcal{H} \hspace{1.3 cm}\label{eqHadSpecHad5B}\\
\mathcal{SH}_k \subset \mathcal{A}_{p,q} \subset \mathcal{A}\label{eqHadSpecHad5C}\\
\mathcal{H}_k \subset \mathcal{A}_{1,q}\hspace{1 cm} \label{eqHadSpecHad5D}
\end{eqnarray}
for all $k>0$, $p \in \left[1,2 \right]$ and $q \in \left[p, p^{\prime} \right]$.

We would like to address the following problems:
\begin{problem}\label{ProblemA1}
We know that linear canonical transforms (LCT) (and in particular the Fourier transform) belong to $\mathcal{SH}_k$ for some $k$ (depending on the particular LCT), and hence also to $\mathcal{H}_k$. Can we devise ways to obtain more operators in $\mathcal{SH}$ and $\mathcal{H}$?
\end{problem}

\vspace{0.1 cm}
\begin{problem}\label{ProblemA2}
Are the inclusions (\ref{eqHadSpecHad5A}-\ref{eqHadSpecHad5D}) proper?
\end{problem}

Let us start by considering Problem \ref{ProblemA1}. Clearly, if $A \in \mathcal{H}$ (or $\mathcal{SH}$), then $\alpha A \in \mathcal{H}$ (or $\mathcal{SH}$), for all $\alpha >0$. It then seems natural to inquire whether these sets are closed with respect to convexity or composition. Unfortunately, the answer is negative.

\begin{lemma}\label{Lemma1}
Let $A,B \in \mathcal{H}$ (or $\mathcal{SH}$). Then, in general $A+B \notin \mathcal{H}$ (or $\mathcal{SH}$).
\end{lemma}

\begin{lemma}\label{Lemma2}
Let $A,B \in \mathcal{H}$ (or $\mathcal{SH}$). Then, in general $AB \notin \mathcal{H}$ (or $\mathcal{SH}$).
\end{lemma}

On a more positive note we have:
\begin{lemma}\label{Lemma3}
Suppose that $U$ and $U^{\ast}$ are bounded operators $L^p(\mathbb{R}^n) \to L^p(\mathbb{R}^n) $ for all $p \in  \left[1, \infty \right]$. Moreover, suppose that
\begin{equation}
U^{\ast}U= \alpha I_{L^p}~,
\label{eqHadSpecHad9}
\end{equation}
for some $\alpha >0$ and all $p \in  \left[1, \infty \right]$. If $A \in \mathcal{H}_k$, $B \in \mathcal{SH}_k$, then $UA \in \mathcal{H}_{\alpha k}$, $UB \in \mathcal{SH}_{\alpha k}$.
\end{lemma}

Here are some examples in connection to Lemma \ref{Lemma3}.

\begin{example}[Multiplication by a phase]\label{ExampleMultiplicationPhase}
Let $\phi: \mathbb{R}^n \to \mathbb{C}$ be a measurable function such that $|\phi(x)|=1$ for almost all $x \in \mathbb{R}^n$. We define the operator:
\begin{equation}
(M_{\phi} f) (x)= \phi(x) f(x)~.
\label{eqHadSpecHad10}
\end{equation}
Then $M_{\phi}$ satisfies the conditions for the operator $U$ (Lemma \ref{Lemma3}) for $\alpha =1$.

\end{example}

\begin{example}[Diffeomorphism]\label{ExampleDiff}

Let $\psi : \mathbb{R}^n \to \mathbb{R}^n$ be a $C^1$-diffeomorphim, with a constant Jacobian:
\begin{equation}
\left| \det \left( \frac{\partial \phi_i}{\partial x_j} \right)_{i,j=1, \cdots, n} \right| =C>0~.
\label{eqHadSpecHad11}
\end{equation}
Define:
\begin{equation}
(D_{\psi} f)(x)=(f \circ \psi )(x)= f \left(\psi(x) \right)~.
\label{eqHadSpecHad12}
\end{equation}
Then $D_{\psi}$ satisfies the conditions for the operator $U$ (Lemma \ref{Lemma3}) for $\alpha =\frac{1}{C}$.

For example, the operator
\begin{equation}
f \in \mathcal{S}(\mathbb{R}^2) \mapsto \int_{\mathbb{R}} \int_{\mathbb{R}} f(x_1,x_2) \exp \left(-2 i \pi x_1 \xi_1 e^{\gamma \xi_1 \xi_2} -2 i \pi x_2 \xi_2 e^{-\gamma \xi_1 \xi_2} \right) dx_1 dx_2~,
\label{eqHadSpecHad13}
\end{equation}
where $\gamma \neq 0$, is a special $1$-Hadamard operator. The reason is that this operator is the composition $D_{\psi} \mathcal{F}$ of the Fourier transform $\mathcal{F}$ (a special $1$-Hadamard operator) and the operator (\ref{eqHadSpecHad12}), where $\psi$ is the diffeomorphism
\begin{equation}
\psi(\xi_1,\xi_2)=(\xi_1 e^{\gamma \xi_1 \xi_2},\xi_2 e^{-\gamma \xi_1 \xi_2})~,
\label{eqHadSpecHad14}
\end{equation}
which has constant Jacobian equal to $1$.
\end{example}

\begin{example}[Composition with step functions]\label{ExampleCompStepFunc}
Let $\left\{\Omega_j,~j \in \mathbb{N}\right\}$ be a (possibly infinite) partition of $\mathbb{R}^n$, that is $\Omega_j$ are measurable subsets of $\mathbb{R}^n$, $\Omega_j \cap \Omega_j = \emptyset$, if $i \neq j$, and $\cup_{j\in \mathbb{N}} \Omega_j = \mathbb{R}^n$.

Let $(\alpha_j)_j$ be a sequence of positive real numbers, such that $0< m \leq \alpha_j \leq M$, $\forall j \in  \mathbb{N}$. We denote by $P_j$ the projection onto $\Omega_j$. If $\chi_j$ is the characteristic function of $\Omega_j$,
\begin{equation}
\chi_j(x)= \left\{
\begin{array}{l l}
1, & \text{if }x \in \Omega_j\\
0, & \text{otherwise}
\end{array}
\right.~, 
\label{eqHadSpecHad15A}
\end{equation}
then $(P_jf)(x)=\chi_j(x) f(x)$. We have of course:
\begin{equation}
P_j^{\ast}=P_j~, ~P_iP_j= \delta_{i,j}P_j~, ~ \forall i,j \in \mathbb{N}~,~ ~~  \sum_{j \in \mathbb{N}} P_j = I_{\mathbb{R}^n}~.
\label{eqHadSpecHad16A}
\end{equation}
Let $A \in \mathcal{SH}_k$, and define:
\begin{equation}
B=\sum_{j\in \mathbb{N}}\alpha_j P_j~.
\label{eqHadSpecHad17A}
\end{equation}
We then have $AB \in \mathcal{H}_{km^2}$. A proof of this can be found in section \ref{ProofExampleStepFunc}.

\end{example}

Concerning Problem \ref{ProblemA2}, we have the following results:

\begin{proposition}\label{PropositionInclusion1}
The inclusions (\ref{eqHadSpecHad5A}) and (\ref{eqHadSpecHad5B}) are proper.
\end{proposition}

\begin{proposition}\label{PropositionInclusion2}
The inclusion (\ref{eqHadSpecHad5C}) is proper for all $k>0$ and all $p \in \left[1,2 \right]$, $q \in \left[p,p^{\prime}\right]$.
\end{proposition}

\subsection{Some open problems}

There remain many open problems to address in this framework of $k$-Hadamard operators. Here is a short list.

\begin{enumerate}
\item In the higher dimensional norm uncertainty principle (Theorem \ref{Theorem2}), we were unable to determine what happens in the threshold cases $q= \frac{2n}{n+2}$ or $q= \frac{2n}{n-2}$, when $n \geq 3$, and $q=1$ or $q= \infty$, when $n=2$. This remains an open problem.

\item In Proposition \ref{Proposition3}, we proved that the higher dimensional norm uncertainty principle cannot hold for $1 \leq q < \frac{2n}{n+2}$ or for $q > \frac{2n}{n-2}$. However, we only proved this if $A= \mathcal{F}$ is the Fourier transform. It would be interesting to determine whether the same result is valid for all special $k$-Hadamard operators.

\item Are those theorems which we proved to be valid for special $k$-Hadamard operators in fact valid for all $k$-Hamamard operators?

\item In Proposition \ref{Proposition2} we were unable to determine whether (\ref{eq8}) holds for $n>2$ and $\frac{2n}{n+2} < q \leq 2$. This case remains an open problem.

\item Can we prove the assumption of Theorem \ref{TheoremEntropicUP} that the constant $C_{p,2,n}$ is differentiable and strictly decreasing near $2$ for all special $k$-Hadamard operators? We know that this claim is valid for the Fourier transform and other metaplectic operators.  
 
\item We believe that Shapiro's uncertainty principle for orthonormal sequences \cite{Jaming,Malinnikova,Shapiro} can be proven with the techniques developed here. Let us first recall the content of Shapiro's theorem. Let $f \in L^2 (\mathbb{R})$ be such that $\|f\|_2=1$. 
In analogy with quantum mechanics, we define the expectation value and the dispersion of $f$ according to:
\[
\begin{array}{l}
\mu (f)= \int_{\mathbb{R}} x |f(x)|^2 dx\\
\\
\Delta (f) = \int_{\mathbb{R}} \left(x- \mu(f) \right)^2 |f(x)|^2 dx
\end{array}
\]
Notice the difference between $V(f)$ and $\Delta(f)$. The two are related by:
\[
V(f)=\Delta(f) + \left(\mu(f) \right)^2 ~.
\]
Then Shapiro's \textit{mean-dispersion principle} states that:

\vspace{0.3 cm}
\textit{There does not exist an infinite orthonormal sequence} $\left\{e_n \right\}_{n=0}^{\infty} \subset L^2 (\mathbb{R})$ \textit{such that all four of }$\mu(e_n)$, $\mu(\widehat{e}_n)$, $\Delta (e_n)$, $\Delta(\widehat{e}_n)$ \textit{are uniformly bounded}.

\vspace{0.3 cm}
We then have from Theorem \ref{Theorem2}:
\[
\begin{array}{c}
\Delta (e_n) + \Delta(\widehat{e}_n) + \mu(e_n) +\mu(\widehat{e}_n)= \\
\\
=V(e_n)+ V(\widehat{e}_n) \geq 2 \sqrt{ V(e_n)V(\widehat{e}_n)} \gtrsim \|e_n\|_q ~\|\widehat{e}_n\|_q
\end{array}
\]
If we can show that $\|e_n\|_q ~ \|\widehat{e}_n\|_q \to \infty$ as $n \to \infty$ for all infinite orthonormal sequences (with some appropriate choice of $q$), then this proves Shapiro's result. Notice that this is true for the Hermite functions for all $1 \leq q < 2$ (see Example \ref{Example1} below).

\end{enumerate}
 
\section*{Acknowledgements}

 The authors would like to thank Y. Tang for some fruitful discussions.

\section*{Notation}

Given $f \in L^1(\mathbb{R}^n)\cap L^2 (\mathbb{R}^n)$, the Fourier transform is defined by:
\begin{equation}
\left(\mathcal{F}f\right)(\xi)= \widehat{f}(\xi)= \int_{\mathbb{R}^n} f(x) e^{-2i\pi x \cdot \xi} dx~,
\label{eqNotation1}
\end{equation}
and extends to $L^2 (\mathbb{R}^n)$. We denote by $\mathcal{S} (\mathbb{R}^n)$ the Schwartz space of test functions, and its dual are the tempered distributions $\mathcal{S}^{\prime} (\mathbb{R}^n)$. The norms $\|f\|_p$ for $p \in \left[\right.1, \infty \left.\right)$ are given by
\begin{equation}
\|f\|_p = \left(\int_{\mathbb{R}^n} |f(x)|^p dx \right)^{\frac{1}{p}}~.
\label{eqNotation2}
\end{equation}
The dimension $n$ is usually clear from the context. Notice that in Theorems \ref{PropositionProof1},\ref{TheoremAdcional2} and \ref{TheoremCowlingPriceInfinity} we admit $p \in \left(0,1 \right)$. We should bear in mind that, in that case, $\|\cdot\|_p$ is not a norm. If $p = \infty$, then:
\begin{equation}
\|f\|_{\infty} = \text{ess sup}_{x \in \mathbb{R}^n} |f(x)|~.
\label{eqNotation3}
\end{equation}
The tensor product of two functions $f,g$ is defined by $(f\otimes g)(x,y)= f(x) g(y)$ and for $n$ functions: $\left(f_1 \otimes \cdots \otimes f_n\right) (x_1,\cdots, x_n)= f(x_1) \cdots f(x_n)$.

\section{Proofs}

In this section, we prove the various results stated in the previous Section.

Let $f $ be some measurable function on $\mathbb{R}^n$, and consider the dilation $f_{\lambda}(x)= f(\lambda x)$, for $\lambda>0$. We have that $f \in L^q(\mathbb{R}^n)$ for some $q \in \left[1, \infty\right]$ if and only if $f_{\lambda} \in L^q(\mathbb{R}^n)$. Moreover, we have:
\begin{equation}
\|f_{\lambda}\|_q= \lambda^{-\frac{n}{q}} \|f\|_q~, ~q \in \left[\right.1,\infty \left. \right)~,
\label{eqExtra1}
\end{equation} 
and
\begin{equation}
\|f_{\lambda}\|_{\infty}= \|f\|_{\infty}~.
\label{eqExtra2}
\end{equation} 
Concerning the variance, we obtain:
\begin{equation}
V(f_{\lambda})=\lambda^{-n-2}V(f)~.
\label{eqExtra3}
\end{equation} 

\subsection{Proof of Proposition \ref{PropositionHausdorffYoung}}

Notice that $A^{\ast} A = k I$, where $I$ is the identity operator on $L^2 (\mathbb{R}^n)$. Thus, $A$ is a $k$-Hadamard operator. We also conclude that $\|A\|_{2 \to 2}= \sqrt{k}$. Since we also have $\|A\|_{1 \to \infty} \leq 1$, it follows from the Riesz-Thorin Theorem that there exists a constant $C^{\prime}>0$, such that:
\begin{equation}
\|Af\|_{p^{\prime}} \leq C^{\prime} \|f\|_p~.
\label{eqProofLemmaHausdorffYoung0}
\end{equation} 

Since $A^{\ast}$ is also a bounded operator $L^1 \to L^{\infty}$ and $L^2 \to L^2$, with the same arguments we conclude that $\|A^{\ast} \|_{p \to p^{\prime}} \leq C^{\prime \prime}$ for some constant $C^{\prime \prime}>0$. 

Since $A f \in L^p (\mathbb{R}^n)$, we have:
\begin{equation}
\|f\|_{p^{\prime}}=\|A^{-1}A f \|_{p^{\prime}}=\frac{1}{k}\|A^{\ast} A f \|_{p^{\prime}} \leq \frac{C^{\prime \prime}}{k}  \|A f \|_{p}~.
\label{eqProofLemmaHausdorffYoung1}
\end{equation} 
If we set $C = \text{max} \left\{C^{\prime},\frac{C^{\prime \prime}}{k} \right\}$, we obtain (\ref{HausdorffYoung2}) and (\ref{HausdorffYoung3}).

\subsection{Proof of Theorem \ref{TheoremAdcional1}}

The result is obviously true if $p=2$, or if  $p<2$ and $p=q$.

We start by remarking that, if $f,Af \in L^p (\mathbb{R}^n)$, then by Hausdorff-Young inequality for special $k$-Hadamard operators (Proposition \ref{PropositionHausdorffYoung}), $f,Af \in L^{p^{\prime}} (\mathbb{R}^n)$:
\begin{equation}
\|f\|_p \gtrsim \|Af\|_{p^{\prime}}, ~~\|Af\|_p \gtrsim \|f\|_{p^{\prime}}.
\label{eqProofTheoremAdicional1.0}
\end{equation}
Hence $f,A f \in L^q (\mathbb{R}^n)$, for all $q \in \left[p,p^{\prime}\right]$, and moreover:
\begin{equation}
\|f\|_p \|A f\|_p \gtrsim  \|f\|_{p^{\prime}} \|Af\|_{p^{\prime}}~. 
\label{eqProofTheoremAdicional1.7}
\end{equation}
This shows that the result is also true if $q=p^{\prime}$. We therefore assume that $p \neq 2$ and $p <q < p^{\prime}$.

Next define 
\begin{equation}
\theta=\frac{p^{\prime}-p}{p^{\prime}-q}~.
\label{eqProofTheoremAdicional1.1}
\end{equation}
Notice that
\begin{equation}
\theta >1, ~~ q-\frac{p}{\theta}=\frac{p^{\prime}(q-p)}{p^{\prime}-p}>0, ~~ \frac{q \theta-p}{\theta-1}=\theta^{\prime} \left(q-\frac{p}{\theta}\right)=p^{\prime}~.
\label{eqProofTheoremAdicional1.2}
\end{equation}
It follows, from H\"older's inequality:
\begin{equation}
\begin{array}{c}
\|f\|_q^q= \int |f|^q = \int |f|^{\frac{p}{\theta}} |f|^{q-\frac{p}{\theta}} \leq \|~ |f|^{\frac{p}{\theta}} \|_{\theta} ~\|~ |f|^{q-\frac{p}{\theta}} \|_{\theta^{\prime}}= \\
\\
=\|f\|_p^{\frac{p}{\theta}}~ \|f\|_{\theta^{\prime} \left(q-\frac{p}{\theta}\right)}^{q-\frac{p}{\theta}} 
=\|f\|_p^{\frac{p}{\theta}}~ \|f\|_{p^{\prime}}^{q-\frac{p}{\theta}}
\end{array}
\label{eqProofTheoremAdicional1.3}
\end{equation}
Upon multiplication by $\|f \|_{p}^{q-\frac{p}{\theta}}$, we obtain:
\begin{equation}
\begin{array}{c}
\|f\|_{p}^{q-\frac{p}{\theta}}~ \|f\|_q^q \leq \|f\|_{p}^{q}~\|f\|_{p^{\prime}}^{q-\frac{p}{\theta}}\\
\\
\Leftrightarrow \left(\frac{\|f\|_p}{\|f\|_{p^{\prime}}}\right)^{q-\frac{p}{\theta}} \leq \left(\frac{\|f\|_p}{\|f\|_q} \right)^{q}  \Leftrightarrow \frac{\|f\|_p}{\|f\|_q} \geq \left(\frac{\|f\|_p}{\|f\|_{p^{\prime}}} \right)^{1-\frac{p}{\theta q}} ~. 
\end{array}
\label{eqProofTheoremAdicional1.4}
\end{equation}
Since
\begin{equation}
1- \frac{p}{\theta q}=\frac{p^{\prime}(q-p)}{q(p^{\prime}-p)}=\frac{1/p-1/q}{1/p-1/p^{\prime}}~,
\label{eqProofTheoremAdicional1.5}
\end{equation}
we obtain:
\begin{equation}
\frac{\|f\|_p}{\|f\|_q} \geq \left(\frac{\|f\|_p}{\|f\|_{p^{\prime}}} \right)^{\frac{1/p-1/q}{1/p-1/p^{\prime}}} ~. 
\label{eqProofTheoremAdicional1.6}
\end{equation}
Finally, from (\ref{eqProofTheoremAdicional1.6}) for $f$ and $A f$, and (\ref{eqProofTheoremAdicional1.7}) the result follows.

Before we conclude, let us make a remark for future reference. Suppose that $A= \mathcal{F}$ is the Fourier transform, and consider the sharp Hausdorff-Young inequality
\begin{equation}
\|\widehat{f} \|_{p^{\prime}} \leq (C_p)^n \|f\|_p ~, ~ \|f \|_{p^{\prime}} \leq (C_p)^n \|\widehat{f}\|_p ~,
\label{eqProofTheoremAdicional1.7A}
\end{equation}
instead of (\ref{eqProofTheoremAdicional1.0}). Here 
\begin{equation}
C_p= \sqrt{\frac{p^{1/p}}{(p^{\prime})^{1/p^{\prime}}}}
\label{eqProofTheoremAdicional1.8}
\end{equation}
is the the Babenko-Beckner constant \cite{Babenko,Beckner,Lieb3}.

Then, following the same steps as before, we conclude, for $p \in \left[1,2 \right]$, $q \in \left[p,p^{\prime} \right]$, that:
\begin{equation}
\mathcal{F}_{p,q}^{(n)} (f) =\frac{\|f\|_p \|\widehat{f}\|_p}{\|f\|_q \|\widehat{f}\|_q}
\geq \left(\frac{1}{C_p}\right)^{\frac{n \left(1/p-1/q\right)}{1/p-1/2}}~, 
\label{eqProofTheoremAdicional1.9}
\end{equation}
for all $f$ with $f,\widehat{f} \in L^p (\mathbb{R}^n)$.

\subsection{Proof of Theorem \ref{TheoremEmbeddings}}

If $p=p^{\prime}=2$, then we have $q=p$ and we have an equality in (\ref{eqEmbeddings3}). So, we assume that $p \neq 2$ and $q \neq p$.

From inequality (\ref{eqProofTheoremAdicional1.6}) and the Hausdorff-Young inequality for special $k$-Hadamard operators (\ref{HausdorffYoung3}), we have:
\begin{equation}
\frac{\|f\|_q}{\|f\|_p} \lesssim \left(\frac{\|f\|_{p^{\prime}}}{\|f\|_p} \right)^{\alpha} \lesssim \left(\frac{\|Af\|_{p}}{\|f\|_p} \right)^{\alpha} ~,
\label{eqProofEmbeddings1}
\end{equation}
where 
\begin{equation}
\alpha = \frac{p^{\prime}(q-p)}{q(p^{\prime}-p)} >0 ~.
\label{eqProofEmbeddings2}
\end{equation}
Let us also remark that:
\begin{equation}
1-\alpha = \frac{p(p^{\prime}-q)}{q(p^{\prime}-p)} \geq 0 ~.
\label{eqProofEmbeddings3}
\end{equation}
If we interchange $Af \longleftrightarrow f$, we conclude also that:
\begin{equation}
\frac{\|Af\|_q}{\|Af\|_p} \lesssim \left(\frac{\|Af\|_{p^{\prime}}}{\|Af\|_p} \right)^{\alpha} \lesssim \left(\frac{\|f\|_{p}}{\|Af\|_p} \right)^{\alpha} ~.
\label{eqProofEmbeddings4}
\end{equation}
From (\ref{eqProofEmbeddings1}), (\ref{eqProofEmbeddings4}), we have:
\begin{equation}
\left\{
\begin{array}{l}
\|f\|_q \lesssim \|Af\|_p^{\alpha} \|f\|_p^{1- \alpha} \\
\\
\|Af\|_q \lesssim \|f\|_p^{\alpha} \|Af\|_p^{1- \alpha}
\end{array}
\right.
\label{eqProofEmbeddings5}
\end{equation}
From the inequality $X^{\alpha}Y^{1-\alpha} \leq \alpha X +(1-\alpha)Y$, which holds for $X,Y>0$ and $0 \leq \alpha \leq 1$, we obtain:
\begin{equation}
\|f\|_q \lesssim M \left(\| f\|_p + \|Af\|_p \right) ~,
\label{eqProofEmbeddings56}
\end{equation}
where $M = \text{max} \left\{\alpha, 1- \alpha\right\}$. 

Using the same arguments for $\|Af\|_q$, we conclude that:
\begin{equation}
\|f\|_q +\|Af\|_q  \lesssim 2M \left( \|f\|_p +\|Af\|_p\right)~. 
\label{eqProofEmbeddings56}
\end{equation}

\subsection{Proof of Theorem \ref{PropositionProof1}}

If $\| |x|^{\theta} f \|_p = \infty$ there is nothing to prove, so we assume that it is finite. Let
\begin{equation}
\widetilde{f}=\frac{f}{\|f\|_{\infty}}~,
\label{eqPropositionProof13}
\end{equation}
and
\begin{equation}
T=\frac{1}{\omega_1^{1/n}} \left(\frac{1}{2}\right)^{\frac{t}{n(t-1)}} \left(\frac{\|f\|_r}{\|f\|_{rt}}\right)^{\frac{rt}{n(t-1)}}~,
\label{eqPropositionProof14}
\end{equation}
with the understanding that
\begin{equation}
T=\left(\frac{1}{2 \omega_1}\right)^{\frac{1}{n}}~ \left(\frac{\|f\|_r}{\|f\|_{\infty}}\right)^{\frac{r}{n}}~,
\label{eqPropositionProof15}
\end{equation}
if $t=\infty$. Here $\omega_1=\frac{2\pi^{n/2}}{n\Gamma (n/2)}$ denotes the volume of the $n$-ball of radius $1$.

Let $B_T= \left\{x \in \mathbb{R}^n:~ |x| \leq T\right\}$ denote the ball of radius $T$. We have, by H\"older's inequality:
\begin{equation}
\begin{array}{c}
\int_{|x|\leq T} |\widetilde{f}|^r = \int \chi_{B_T} |\widetilde{f}|^r\leq \|\chi_{B_T}\|_{\frac{t}{t-1}} \||\widetilde{f}|^r\|_t =\\
\\
=(\omega_1 T^n)^{\frac{t-1}{t}} \|\widetilde{f}\|_{tr}^r = \frac{1}{2}\|\widetilde{f}\|_r^r~.
\end{array}
\label{eqPropositionProof16}
\end{equation}
In view of the previous inequality, H\"older's inequality and the fact that $|\widetilde{f}|\leq 1$ a.e., we have for all $p,u,\theta>0$ satisfying (\ref{eqPropositionProof11}) that:
\begin{equation}
\begin{array}{c}
\frac{1}{2}\|\widetilde{f}\|_r^r\leq \int_{|x|>T} |\widetilde{f}|^r \leq \int_{|x|>T} |\widetilde{f}|^u =\\
\\
= \int_{|x|>T}\frac{1}{|x|^{\theta u}}\left(|x|^{\theta u}|\widetilde{f}|^u\right)\\
\\
 \leq \left(\int_{|x|>T}\frac{1}{|x|^{\frac{\theta up}{p-u}}} dx\right)^{\frac{p-u}{p}} 
\left(\int_{\mathbb{R}^n}|x|^{\theta p} |\widetilde{f}(x)|^p dx\right)^{\frac{u}{p}} =\\
\\
=\left(\sigma_1 \int_T^{\infty} r^{n-1-\frac{\theta up}{p-u}} dr \right)^{\frac{p-u}{p}} \| |x|^{\theta} \widetilde{f} \|_p^u~, 
\end{array}
\label{eqPropositionProof17}
\end{equation}
where $\sigma_1=\frac{2\pi^{n/2}}{\Gamma (n/2)}$ is the surface area of the $n$-sphere of radius $1$. The integral on the right-hand side of (\ref{eqPropositionProof17}) is convergent, provided $\theta>n \left(\frac{1}{u}-\frac{1}{p}\right)$, in which case we have:
\begin{equation}
\int_T^{\infty} r^{n-1-\frac{\theta up}{p-u}} dr = \frac{p-u}{\theta up-n(p-u)} T^{\frac{n(p-u)-\theta up}{p-u}}~.
\label{eqPropositionProof17A}
\end{equation}
It follows that:
\begin{equation}
\begin{array}{c}
\||x|^{\theta}  \widetilde{f} \|_p\gtrsim \|\widetilde{f} \|_r^{r/u} T^{\frac{u\theta p-n(p-u)}{up}} \\
\\
\Longleftrightarrow 
\||x|^{\theta}  \widetilde{f} \|_p\gtrsim \|\widetilde{f} \|_r^{r/u} \left(\frac{\|f\|_r}{\|f\|_{rt}}\right)^{\frac{rt(u \theta p-n(p-u))}{(t-1)npu}}\\
\\
\Longleftrightarrow \||x|^{\theta}  f \|_p\gtrsim \|f\|_{\infty} \left(\frac{\|f \|_r}{\|f\|_{\infty}}\right)^{r/u} \left(\frac{\|f\|_r}{\|f\|_{rt}}\right)^{\frac{rt(u \theta p-n(p-u))}{(t-1)npu}}
\end{array}
\label{eqPropositionProof18}
\end{equation}
Notice that all norms appearing in the proof are finite, since by assumption, $f \in L^r \cap L^{\infty}$ when $u \neq r$, and $f \in L^r \cap L^{rt}$ when $u=r$.

The proof for $u=r$ follows the same steps as before, but we do not need to consider $\widetilde{f}$ and hence we do not require $f \in L^{\infty} $, but only $f \in L^r \cap L^{rt}$. 

\subsection{Proof of Theorem \ref{TheoremAdcional2}}

As previously, let
\begin{equation}
T= \frac{1}{\omega_1^{\frac{1}{n}}} \left(\frac{1}{2}\right)^{\frac{t}{n(t-1)}} \left(\frac{\|f\|_p}{\|f\|_{pt}}\right)^{\frac{tp}{n(t-1)}}~,
\label{eqProofTheoremAdcional2.1}
\end{equation}
if $1< t < \infty$, and
\begin{equation}
T=\left(\frac{1}{2 \omega_1} \right)^{\frac{1}{n}} \left(\frac{\|f\|_p}{\|f\|_{\infty}}\right)^{\frac{p}{n}}~,
\label{eqProofTheoremAdcional2.2}
\end{equation}
if $t = \infty$.

Following the same steps as in the proof of the previous theorem, we conclude that:
\begin{equation}
\int_{|x|\leq T} |f|^p \leq \frac{1}{2} \|f\|_p^p ~.
\label{eqProofTheoremAdcional2.3}
\end{equation}
Consequently:
\begin{equation}
\frac{1}{2} \|f\|_p^p \leq \int_{|x| > T} |f|^p =\int_{|x| > T}\frac{1}{|x|^{\theta p}} \left(|x|^{\theta p} |f|^p \right) \leq \frac{1}{T^{\theta p}} \|~|x|^{\theta} f \|_p^p ~.
\label{eqProofTheoremAdcional2.4}
\end{equation}
Equivalently:
\begin{equation}
\|~|x|^{\theta} f \|_p \gtrsim T^{\theta} \|f\|_p \Leftrightarrow \|~|x|^{\theta} f \|_p \gtrsim \left(\frac{\|f\|_p}{\|f\|_{tp}}\right)^{\frac{\theta p t}{n(t-1)}} \|f\|_p~.
\label{eqProofTheoremAdcional2.5}
\end{equation}
If we set $q=tp$, with $q \in \left( \right. p, \infty \left.\right]$, the result follows.

\subsection{Proof of Corollary \ref{CorollaryAdicional1}}

If $\|\left(- \Delta \right)^s f \|_2 = \infty $ there is nothing to prove. So suppose that$\|\left(- \Delta \right)^s f \|_2 =\frac{1}{2\pi}\|~| \xi|^{2s} \widehat{f} \|_2 < \infty $, and let $1 \leq p<2$. 

From Theorem \ref{TheoremAdcional2}, we have (setting $q=p^{\prime}$ in (\ref{eqTheoremAdcional2.1})):
\begin{equation}
 \|\left(- \Delta \right)^s f \|_2  \gtrsim \left(\frac{\|f\|_2}{\|\widehat{f}\|_{p^{\prime}}} \right)^{\frac{4sp^{\prime}}{n(p^{\prime}-2)}} \|f\|_2~.
\label{eqProofCorollaryAdicional1.3}
\end{equation}
From the Hausdorff-Young inequality:
\begin{equation}
 \|\left(- \Delta \right)^s f \|_2  \gtrsim \left(\frac{\|f\|_2}{\|f\|_{p}} \right)^{\frac{4sp}{n(2-p)}} \|f\|_2~ \Leftrightarrow  \|\left(- \Delta \right)^s f \|_2  \gtrsim \left(\frac{\|f\|_p}{\|f\|_2} \right)^{\frac{1}{\alpha}} \|f\|_2~.
\label{eqProofCorollaryAdicional1.4}
\end{equation}

\subsection{Proof of Theorem \ref{TheoremCowlingPriceInfinity}}

In (\ref{eqPropositionProof12}) set $r=1$, $t= \infty$, so that
\begin{equation}
0<u \leq 1 ~, ~p>u~, ~\theta > n \left(\frac{1}{u}-\frac{1}{p}\right)~.
\label{eqProofTheoremCowlingPriceInfinity1}
\end{equation}
We obtain:
\begin{equation}
\||x|^{\theta}f \|_p \gtrsim\|f\|_{\infty} \left(\frac{\|f\|_1}{\|f\|_{\infty}}\right)^{\frac{\theta}{n} + \frac{1}{p}}~.
\label{eqProofTheoremCowlingPriceInfinity2}
\end{equation}
The analogous inequality for $Af$, with 
\begin{equation}
0<v \leq 1 ~, ~q>v~, ~\phi > n\left(\frac{1}{v}-\frac{1}{q}\right)~,
\label{eqProofTheoremCowlingPriceInfinity3}
\end{equation}
is
\begin{equation}
\||\xi|^{\phi}Af \|_q \gtrsim\|Af\|_{\infty} \left(\frac{\|Af\|_1}{\|Af\|_{\infty}}\right)^{\frac{\phi}{n} + \frac{1}{q}}~.
\label{eqProofTheoremCowlingPriceInfinity4}
\end{equation}
Given $\alpha, \beta >0$, we obtain from (\ref{eqProofTheoremCowlingPriceInfinity2},\ref{eqProofTheoremCowlingPriceInfinity4}):
\begin{equation}
\||x|^{\theta}f \|_p^{\alpha}~ \||\xi|^{\phi}Af \|_q^{\beta} \gtrsim\|f\|_{\infty}^{\alpha}\|Af\|_{\infty}^{\beta} \left(\frac{\|f\|_1}{\|f\|_{\infty}}\right)^{\alpha \left(\frac{\theta}{n} + \frac{1}{p}\right)} \left(\frac{\|Af\|_1}{\|Af\|_{\infty}}\right)^{\beta \left(\frac{\phi}{n} + \frac{1}{q}\right)}~.
\label{eqProofTheoremCowlingPriceInfinity5}
\end{equation}
On the other hand, if $\alpha \left(\frac{\theta}{n} + \frac{1}{p}\right)=\beta \left(\frac{\phi}{n} + \frac{1}{q}\right)$, then:
\begin{equation}
\||x|^{\theta}f \|_p^{\alpha}~ \||\xi|^{\phi}Af \|_q^{\beta} \gtrsim\|f\|_{\infty}^{\alpha}\|Af\|_{\infty}^{\beta} \left(\frac{\|f\|_1 \|Af\|_1}{\|f\|_{\infty} \|Af\|_{\infty}}\right)^{\alpha \left(\frac{\theta}{n} + \frac{1}{p}\right)}~.
\label{eqProofTheoremCowlingPriceInfinity6}
\end{equation}
From the primary uncertainty principle (\ref{eq0A}), it follows that:
\begin{equation}
\||x|^{\theta}f \|_p^{\alpha}~ \||\xi|^{\phi}Af \|_q^{\beta} \gtrsim\|f\|_{\infty}^{\alpha}\|Af\|_{\infty}^{\beta} ~.
\label{eqProofTheoremCowlingPriceInfinity7}
\end{equation}
Finally, suppose that $\theta,p >0$ are such that $\frac{\theta}{n}+\frac{1}{p}>1$. If $p \geq 1$, we can always find $0<u<1$, such that $\frac{\theta}{n}+\frac{1}{p}> \frac{1}{u} >1 \geq \frac{1}{p}$. 
Alternatively, if $0<p<1$, let $0 < \epsilon < p \left(1- \frac{n}{\theta p+n} \right)$ and set $u=p- \epsilon$. It follows that $0 <u <p <1$ and $\frac{\theta}{n}+\frac{1}{p}> \frac{1}{u}$.  Altogether, this means that under the conditions of the theorem, we can always find $0 <u<1$ such that (\ref{eqProofTheoremCowlingPriceInfinity1}) holds. The same can be said about $0<v<1$.

\subsection{Proof of Theorem \ref{TheoremCowlingPriceHadamard}}

In (\ref{eqPropositionProof12}) set $u=r=1$, so that
\begin{equation}
1<t \leq \infty ~, ~p>1~, ~\theta >n \left( 1-\frac{1}{p} \right)~.
\label{eqProofTheoremCowlingPriceHadamard1}
\end{equation}
We obtain:
\begin{equation}
\||x|^{\theta}f \|_p \gtrsim\|f\|_1 \left(\frac{\|f\|_1}{\|f\|_t}\right)^{\frac{t(\theta p-n(p-1))}{np(t-1)}}~.
\label{eqProofTheoremCowlingPriceHadamard2}
\end{equation}
The analogous inequality for $Af$, with 
\begin{equation}
q>1~, ~\phi >n \left( 1-\frac{1}{q} \right)~,
\label{eqProofTheoremCowlingPriceHadamard3}
\end{equation}
is
\begin{equation}
\||\xi|^{\phi}Af \|_q \gtrsim\|Af\|_1 \left(\frac{\|Af\|_1}{\|Af\|_t}\right)^{\frac{t(\phi q-n( q-1))}{nq(t-1)}}~.
\label{eqProofTheoremCowlingPriceHadamard4}
\end{equation}
Given $\alpha, \beta >0$, we obtain from (\ref{eqProofTheoremCowlingPriceHadamard2},\ref{eqProofTheoremCowlingPriceHadamard4}):
\begin{equation}
\||x|^{\theta}f \|_p^{\alpha}~ \||\xi|^{\phi}Af \|_q^{\beta} \gtrsim\|f\|_1^{\alpha}\|Af\|_1^{\beta} \left(\frac{\|f\|_1}{\|f\|_t}\right)^{\frac{\alpha t(\theta p-n(p-1))}{np(t-1)}} \left(\frac{\|Af\|_1}{\|Af\|_t}\right)^{\frac{\beta t(\phi q-n(q-1))}{nq(t-1)}}~.
\label{eqProofTheoremCowlingPriceHadamard5}
\end{equation}
On the other hand, if $\alpha \left(\theta -n\left(1- \frac{1}{p}\right) \right)=\beta \left(\phi -n \left(1- \frac{1}{q}\right)\right)$, then:
\begin{equation}
\||x|^{\theta}f \|_p^{\alpha}~ \||\xi|^{\phi}Af \|_q^{\beta} \gtrsim\|f\|_1^{\alpha}\|Af\|_1^{\beta} \left(\frac{\|f\|_1 \|Af\|_1}{\|f\|_t \|Af\|_t}\right)^{\frac{\alpha t(\theta p-n(p-1))}{np(t-1)}}~.
\label{eqProofTheoremCowlingPriceHadamard6}
\end{equation}
From the norm uncertainty principle (\ref{eq0B}), it follows that:
\begin{equation}
\||x|^{\theta}f \|_p^{\alpha}~ \||\xi|^{\phi}Af \|_q^{\beta} \gtrsim\|f\|_1^{\alpha}\|Af\|_1^{\beta} ~.
\label{eqProofTheoremCowlingPriceHadamard7}
\end{equation}
If $\alpha =\beta$, we obtain from (\ref{eq0A},\ref{eq0B}) and (\ref{eqProofTheoremCowlingPriceHadamard7}):
\begin{equation}
\||x|^{\theta}f \|_p~ \||\xi|^{\phi}Af \|_q \gtrsim\|f\|_s \|Af\|_s ~,
\label{eqProofTheoremCowlingPriceHadamard8}
\end{equation} 
for all $s \in \left[1, \infty \right]$.

\subsection{Proof of Theorem \ref{TheoremAdicional3}}

From Theorem \ref{TheoremAdcional2}, we have for $\theta>0$, $1 \leq p <2$ and $p < r \leq p^{\prime}$:
\begin{equation}
\|~|x|^{\theta} f \|_p \gtrsim \left(\frac{\|f\|_p}{\|f\|_r}\right)^{\frac{\theta r p}{n(r-p)}} \|f\|_p~,
\label{eqProofTheoremAdicional3.1}
\end{equation}
 and
\begin{equation}
\|~|\xi|^{\theta} A f \|_p \gtrsim \left(\frac{\|A f\|_p}{\| A f\|_r}\right)^{\frac{\theta r p}{n(r-p)}} \|A f\|_p~.
\label{eqProofTheoremAdicional3.2}
\end{equation}
Multiplying the two inequalities yields:
\begin{equation}
\|~|x|^{\theta} f \|_p\|~|\xi|^{\theta} A f \|_p \gtrsim \left(\frac{\|f\|_p \|A f\|_p}{\|f\|_ r \|A f\|_r}\right)^{\frac{\theta r p}{n(r-p)}} \|f\|_p \|A f\|_p~.
\label{eqProofTheoremAdicional3.3}
\end{equation}
From Theorem \ref{TheoremAdcional1}, we obtain:
\begin{equation}
\|~|x|^{\theta} f \|_p\|~|\xi|^{\theta} A f \|_p \gtrsim  \|f\|_p \|A f\|_p~.
\label{eqProofTheoremAdicional3.4}
\end{equation}
Finally, using Theorem \ref{TheoremAdcional1} again, we have for $q \in \left[p,p^{\prime}\right]$:
\begin{equation}
\|~|x|^{\theta} f \|_p\|~|\xi|^{\theta} A f \|_p \gtrsim  \|f\|_q \|A f\|_q~.
\label{eqProofTheoremAdicional3.5}
\end{equation}

\subsection{Proof of Theorem \ref{TheoremPreCowlingPrice}}

Setting $u=r$ in Theorem \ref{PropositionProof1} we obtain:
\begin{equation}
\||x|^{\theta}f \|_p \gtrsim \|f\|_r \left(\frac{\|f\|_r}{\|f\|_{rt}} \right)^{\frac{t\left(r \theta p-n(p-r)\right)}{np(t-1)}}~,
\label{eqProofTheoremPreCowlingPrice1}
\end{equation}
and likewise:
\begin{equation}
\||\xi|^{\phi} Af \|_q \gtrsim \|Af\|_r \left(\frac{\|Af\|_r}{\|Af\|_{rt}} \right)^{\frac{t\left(r \phi q-n(q-r)\right)}{nq(t-1)}}~.
\label{eqProofTheoremPreCowlingPrice2}
\end{equation}
Thus, for any $\alpha, \beta >0$:
\begin{equation}
\||x|^{\theta}f \|_p^{\alpha} \||\xi|^{\phi} Af \|_q^{\beta} \gtrsim \|f\|_r^{\alpha} \|Af\|_r^{\beta} \left(\frac{\|f\|_r}{\|f\|_{rt}} \right)^{\frac{\alpha t\left(r \theta p-n(p-r)\right)}{np(t-1)}}  \left(\frac{\|Af\|_r}{\|Af\|_{rt}} \right)^{\frac{t \beta \left(r \phi q-n(q-r)\right)}{nq(t-1)}}~.
\label{eqProofTheoremPreCowlingPrice3}
\end{equation}
If
\begin{equation}
\alpha \left[\theta-n \left(\frac{1}{r}-\frac{1}{p}\right) \right]=\beta\left[\phi-n \left(\frac{1}{r}-\frac{1}{q}\right) \right]~,
\label{eqProofTheoremPreCowlingPrice4}
\end{equation}
then:
\begin{equation}
\||x|^{\theta}f \|_p^{\alpha} \||\xi|^{\phi} Af \|_q^{\beta} \gtrsim \|f\|_r^{\alpha} \|Af\|_r^{\beta} \left(\frac{\|f\|_r \|Af\|_r}{\|f\|_{rt} \|Af\|_{rt}} \right)^{\frac{\alpha t\left(r \theta p-n(p-r)\right)}{np(t-1)}}~.
\label{eqProofTheoremPreCowlingPrice5}
\end{equation}
If we now assume that $1 < t \leq r^{\prime}/r$, then we conclude, from Theorem \ref{TheoremAdcional1}, that:
\begin{equation}
\||x|^{\theta}f \|_p^{\alpha} \||\xi|^{\phi} Af \|_q^{\beta} \gtrsim \|f\|_r^{\alpha} \|Af\|_r^{\beta}~.
\label{eqProofTheoremPreCowlingPrice6}
\end{equation}

\subsection{Proof of Theorem \ref{Theorem2}}

In Theorem \ref{TheoremPreCowlingPrice}, set $\theta=\phi=1$, $p=q=2$, $\alpha=\beta$ and $1 \leq r <2$. Then, in view of (\ref{eqPreCowlingPrice1}), we must also have:
\begin{equation}
r > \frac{2n}{n+2}~.
\label{eqProofTheorem21}
\end{equation}
We thus obtain:
\begin{equation}
V(f) V(Af) \gtrsim \|f\|_r^2 \|Af\|_r^2~,
\label{eqProofTheorem22}
\end{equation}
for all $f$, such that $f,Af \in L^r (\mathbb{R}^n)$.

Next, choose $q \in \left[r,r^{\prime} \right]$. From Theorem \ref{TheoremAdcional1}, we conclude that:
\begin{equation}
V(f) V(Af) \gtrsim \|f\|_q^2 \|Af\|_q^2~,
\label{eqProofTheorem23}
\end{equation}
for all $f$, such that $f, Af \in L^r (\mathbb{R}^n)$.

\begin{enumerate}
\item Suppose that $n=1$. Then $\frac{2}{3} < r < 2$. Thus, if we choose $r=1$, $r^{\prime}= \infty$, we conclude that (\ref{eqProofTheorem23}) holds for all $q \in \left[1, \infty\right]$.  

\item Let $n=2$. Then $1 <r<2$, and (\ref{eqProofTheorem23}) holds for all $q \in \left(1, \infty\right)$. 

\item Finally, suppose that $n \geq 3$. Then $\frac{2n}{n+2}<r<2$. Since the H\"older dual of $\frac{2n}{n+2}$ is $\frac{2n}{n-2}$, if follows that (\ref{eqProofTheorem23}) holds for all $q \in \left(\frac{2n}{n+2}, \frac{2n}{n-2} \right)$. 
\end{enumerate}

\subsection{Proof of Proposition \ref{Proposition3}}

Suppose first that $q \in \left[ \right. 1, \frac{2n}{n+2} \left. \right)$. Consider the function $h_{\alpha}(x)= \prod_{j=1}^n h_k(x_j)$, for some $k =0,1,2,3,\cdots$, where $h_k$ is the Hermite function (\ref{eqPrep1}) and $\alpha =(k,k, \cdots,k)$ is a multi-index, with $|\alpha|=nk$. In view of (\ref{eqPrep2},\ref{eqPrep5},\ref{eqPrep6}), we have:
\begin{equation}
V(h_{\alpha})=(2\pi)^2 V(\widehat{h}_{\alpha})= E_{\alpha}=|\alpha|+\frac{n}{2}= n \left(k+\frac{1}{2}\right)~.
\label{eqProofProposition31}
\end{equation}
Since, by assumption $q<2$, we have from the asymptotic behaviour in Example \ref{Example1} and (\ref{eqProofProposition31})
\begin{equation}
\begin{array}{c}
\frac{V(h_{\alpha}) V(\widehat{h}_{\alpha})}{\|h_{\alpha}\|_q^2\|\widehat{h}_{\alpha}\|_q^2}= \frac{1}{(2 \pi)^{2+n-2n/q}} \left(\frac{V(h_{\alpha}) }{\|h_{\alpha}\|_q^2}\right)^2=\\
\\
= \frac{1}{(2 \pi)^{2+n-2n/q}}\left(\frac{n \left(k+\frac{1}{2}\right)}{\|h_{k}\|_q^{2n}}\right)^2 \sim  \frac{1}{(2 \pi)^{2+n-2n/q}} \left(\frac{n \left(k+\frac{1}{2}\right)}{k^{n\left(1/q-1/2\right)}}\right)^2
\end{array}
\label{eqProofProposition32}
\end{equation}
We conclude that, if $q< \frac{2n}{n+2}$, the previous expression goes to zero as $k \to + \infty$.

Alternatively, suppose that $q \in \left(\frac{2n}{n-2}, \infty \right)$, and consider the function $g_c^{(n)}$ of Example \ref{Example2}.

Since $\widehat{g}_c^{(n)}=g_c^{(n)}$, we have:
\begin{equation}
\frac{V(g_c^{(n)}) V(\widehat{g}_c^{(n)})}{\|g_c^{(n)}\|_q^2 \|\widehat{g}_c^{(n)}\|_q^2} =\left(\frac{V(g_c^{(n)}) }{\|g_c^{(n)}\|_q^2}\right)^2 ~.
\label{eqProofProposition33}
\end{equation}
As $c \to 0^+$, we have from (\ref{eqPrep11})

\begin{equation}
\begin{array}{c}
V(g_{c}^{(n)} )=\frac{n}{\pi} \left(\sqrt{2} +\frac{2c}{\sqrt{c^4+1}} \right)^{n-1}\left[\frac{1}{4\sqrt{2}}\left(\frac{c^4+1}{c^2}\right) +\frac{c^3}{\left(c^4+1\right)^{3/2}} \right]\\
\\
\sim \frac{n(\sqrt{2})^{n-6}}{\pi c^2} \left[1+ (n-1) \sqrt{2} c + \mathcal{O} (c^2) \right]~.
\end{array}
\label{eqProofProposition34}
\end{equation}  
Consequently, from (\ref{eqProofProposition33},\ref{eqProofProposition34},\ref{eqPrep12}), we obtain:
\begin{equation}
\begin{array}{c}
\frac{V(g_c^{(n)}) V(\widehat{g}_c^{(n)})}{\|g_c^{(n)}\|_q^2 \|\widehat{g}_c^{(n)}\|_q^2} \leq\\
\\
\leq  \left\{\frac{n 4^n (\sqrt{2})^{n-6} q^{n/q}}{\pi} \frac{\left[1+ (n-1) \sqrt{2} c + \mathcal{O} (c^2) \right]}{(1+c^{1-2/q})^{2n}}  c ^{n(1-2/q)-2} \right\}^2~.
\end{array}
\label{eqProofProposition35}
\end{equation}  
Thus, if $ \infty>q > \frac{2n}{n-2}$, the right-hand side goes to zero, as $c\to 0^+$.

Finally, if $q= \infty$, we have from (\ref{eqProofProposition33},\ref{eqProofProposition34},\ref{eqPrep9})
\begin{equation}
\begin{array}{c}
\frac{V(g_c^{(n)}) V(\widehat{g}_c^{(n)})}{\|g_c^{(n)}\|_{\infty}^2 \|\widehat{g}_c^{(n)}\|_{\infty}^2} \leq\\
\\
\leq   \left\{\frac{n (\sqrt{2})^{n-6} }{\pi} \left[1+ \left((n-1) \sqrt{2} -2n\right)c + \mathcal{O} (c^2) \right] c ^{n-2} \right\}^2~. 
\end{array}
\label{eqProofProposition36}
\end{equation}  
Again, this vanishes as $c\to 0^+$.

\subsection{Proof of Theorem \ref{TheoremEntropicUP}}

Let us start by defining the functionals:
\begin{equation}
I_f (p):= \left(\frac{\|f\|_p \|A f\|_p}{\|f\|_2 \|Af\|_2}\right)^p ~, ~ 1 \leq p \leq 2~.
\label{eqProofEntropic1}
\end{equation}
If $p=2$, then of course $I_f(2)=1$ for all $f$. If $1 \leq p<2$ and $q \in \left[p, p^{\prime} \right]$, we have from Theorem \ref{TheoremAdcional1}:
\begin{equation}
I_f(p)=\left(\frac{\|f\|_p \|A f\|_p}{\|f\|_q \|Af\|_q} \frac{\|f\|_q \|Af\|_q}{\|f\|_2 \|Af\|_2}\right)^p \geq C_{p,q,n}^p \left(\frac{\|f\|_q \|A f\|_q}{\|f\|_2 \|Af\|_2}\right)^p~.
\label{eqProofEntropic1A}
\end{equation}
Consequently, if we set $q=2$:
\begin{equation}
I_f(p)\geq C_{p,2,n}^p ~.
\label{eqProofEntropic1B}
\end{equation}
Since, by assumption, $C_{p,2,n}$ is differentiable and strictly decreasing near $2$ with respect to p, we have:
\begin{equation}
I_f(p)\geq 1+ \frac{C_n}{2} (2-p)+ \mathcal{O} (2-p)^2 ~,
\label{eqProofEntropic1C}
\end{equation}
with $C_n=-2 \left. \frac{\partial C_{p,n,2}}{\partial p}\right|_{p=2}>0$, and where we used the fact that $C_{2,2,n}=1$.

Let $\epsilon>0$, and define
\begin{equation}
K(\epsilon):= \frac{1}{\epsilon} \left[I_f (2 -2 \epsilon) - I_f (2) \right]~.
\label{eqProofEntropic3}
\end{equation}
From (\ref{eqProofEntropic1B},\ref{eqProofEntropic1C},\ref{eqProofEntropic3}), we have:
\begin{equation}
K (\epsilon) \geq \frac{1}{\epsilon} \left[1+C_n \epsilon + \mathcal{O}(\epsilon^2)-1 \right]= C_n + \mathcal{O}(\epsilon)~.
\label{eqProofEntropic4A}
\end{equation}

We can now prove (\ref{eqTheoremEntropicUP2}).

For any $X >0$, we have (see \cite{Lieb1}):
\begin{equation}
1+ \epsilon \ln X \leq X^{\epsilon}~,
\label{eqProofEntropic5}
\end{equation}
which is equivalent to:
\begin{equation}
\frac{1}{\epsilon} X \left(X^{- \epsilon}-1 \right) \leq - X^{1- \epsilon} \ln X~.
\label{eqProofEntropic6}
\end{equation}
Setting 
\begin{equation}
X= \left|\frac{f(x) Af(\xi)}{\|f\|_2 \|Af \|_2} \right|^2
\label{eqProofEntropic7}
\end{equation}
in the previous inequality, we obtain:
\begin{equation}
\frac{1}{\epsilon} \left(\left|\frac{f(x) Af(\xi)}{\|f\|_2 \|Af \|_2} \right|^{2-2 \epsilon}-\left|\frac{f(x) Af(\xi)}{\|f\|_2 \|Af \|_2} \right|^2 \right) \leq - \left|\frac{f(x) Af(\xi)}{\|f\|_2 \|Af \|_2} \right|^{2- 2\epsilon} \ln \left|\frac{f(x) Af(\xi)}{\|f\|_2 \|Af \|_2} \right|^2~.
\label{eqProofEntropic8}
\end{equation}
Upon integration over $x$ and $\xi$, we obtain:
\begin{equation}
\begin{array}{c}
K(\epsilon) \leq - \int_{\mathbb{R}^n} \int_{\mathbb{R}^n} \left|\frac{f(x) Af(\xi)}{\|f\|_2 \|Af \|_2} \right|^{2- 2\epsilon} \ln \left|\frac{f(x) Af(\xi)}{\|f\|_2 \|Af \|_2} \right|^2 dx d \xi\\
\\
\Leftrightarrow K(\epsilon) \leq - \left(\int_{\mathbb{R}^n} \left|\frac{Af(\xi)}{\|Af \|_2} \right|^{2- 2\epsilon} d \xi \right) \cdot \left(  \int_{\mathbb{R}^n} \left|\frac{f(x) }{\|f\|_2 } \right|^{2- 2\epsilon} \ln \left|\frac{f(x) }{\|f\|_2 } \right|^2 dx  \right) +\\
\\
+ \left(f\longleftrightarrow Af \right)~.
\end{array}
\label{eqProofEntropic9}
\end{equation}
Now, let us consider the various terms appearing on the right-hand side of the previous inequality. In the sequel we denote by $\Omega$ and $\widehat{\Omega}$ the sets:
\begin{equation}
\begin{array}{l}
\Omega:= \left\{ x \in \mathbb{R}^n:~\frac{|f(x)|}{\|f\|_2} \geq 1\right\}\\
\\
\widehat{\Omega}:= \left\{ \xi \in \mathbb{R}^n:~\frac{|Af(\xi)|}{\|Af\|_2} \geq 1\right\}
\end{array}
\label{eqProofEntropic10}
\end{equation}
and by $\Omega^c, \widehat{\Omega}^c$ their complements.

For $0 < \epsilon < 1/2$:
\begin{equation}
\begin{array}{c}
\int_{\mathbb{R}^n} \left|\frac{Af(\xi)}{\|Af \|_2} \right|^{2- 2\epsilon} d \xi = \int_{\widehat{\Omega}} \left|\frac{Af(\xi)}{\|Af \|_2} \right|^{2- 2\epsilon} d \xi + \int_{\widehat{\Omega}^c} \left|\frac{Af(\xi)}{\|Af \|_2} \right|^{2- 2\epsilon} d \xi \\
\\
\leq \int_{\widehat{\Omega}} \left|\frac{Af(\xi)}{\|Af \|_2} \right|^{2} d \xi + \int_{\widehat{\Omega}^c} \left|\frac{Af(\xi)}{\|Af \|_2} \right| d \xi \leq 1 + \frac{\|Af\|_1}{\|Af\|_2}~.
\end{array}
\label{eqProofEntropic11}
\end{equation}
Thus, by dominated convergence:
\begin{equation}
\int_{\mathbb{R}^n} \left|\frac{Af(\xi)}{\|Af \|_2} \right|^{2- 2\epsilon} d \xi \to 1~, \text{ as } \epsilon \to 0^+. 
\label{eqProofEntropic12}
\end{equation}
Similarly:
\begin{equation}
\begin{array}{c}
\left|- \int_{\mathbb{R}^n} \left|\frac{f(x) }{\|f\|_2 } \right|^{2- 2\epsilon} \ln \left|\frac{f(x) }{\|f\|_2 } \right|^2 dx \right|\\
\\
\leq \left|- \int_{\Omega} \left|\frac{f(x) }{\|f\|_2 } \right|^{2- 2\epsilon} \ln \left|\frac{f(x) }{\|f\|_2 } \right|^2 dx \right|+\left|- \int_{\Omega^c} \left|\frac{f(x) }{\|f\|_2 } \right|^{2- 2\epsilon} \ln \left|\frac{f(x) }{\|f\|_2 } \right|^2 dx \right|=\\
\\
= \int_{\Omega} \left|\frac{f(x) }{\|f\|_2 } \right|^{2- 2\epsilon} \ln \left|\frac{f(x) }{\|f\|_2 } \right|^2 dx - \int_{\Omega^c} \left|\frac{f(x) }{\|f\|_2 } \right|^{2- 2\epsilon} \ln \left|\frac{f(x) }{\|f\|_2 } \right|^2 dx \\
\\
\leq \int_{\Omega} \left|\frac{f(x) }{\|f\|_2 } \right|^{2} \ln \left|\frac{f(x) }{\|f\|_2 } \right|^2 dx -  \int_{\Omega^c} \left|\frac{f(x) }{\|f\|_2 } \right| \ln \left|\frac{f(x) }{\|f\|_2 } \right|^2 dx < \infty
\end{array}
\label{eqProofEntropic13}
\end{equation}
Again, by dominated convergence:
\begin{equation}
- \int_{\mathbb{R}^n} \left|\frac{f(x) }{\|f\|_2 } \right|^{2- 2\epsilon} \ln \left|\frac{f(x) }{\|f\|_2 } \right|^2 dx \to - \int_{\mathbb{R}^n} \left|\frac{f(x) }{\|f\|_2 } \right|^{2} \ln \left|\frac{f(x) }{\|f\|_2 } \right|^2 dx = H \left[f\right] ~,
\label{eqProofEntropic14}
\end{equation}
as $\epsilon \to 0^+$.

Similar results hold, when we interchange $f \longleftrightarrow Af$. Altogether, the result follows from (\ref{eqProofEntropic4A},\ref{eqProofEntropic9},\ref{eqProofEntropic12},\ref{eqProofEntropic14}), as $\epsilon \to 0^+$.

\subsection{Proof of Proposition \ref{Proposition1}}

Let us consider the functional 
\begin{equation}
\mathcal{G}_{\beta_n,\infty}^{(n)} (f)= \left(\frac{V(f)}{\|f\|_{\infty}^2}\right)^{\beta_n} \frac{\|f\|_{\infty}}{\|f\|_1}~,
\label{eqProposition21}
\end{equation}
for $\beta_n>0$ and $f \in \mathcal{S}(\mathbb{R}^n)$. 

Using the dilation $f_{\lambda}(x)=f(\lambda x)$, we obtain from (\ref{eqExtra1}-\ref{eqExtra3}):
\begin{equation}
\mathcal{G}_{\beta_n,\infty}^{(n)} (f_{\lambda} )=\lambda^{-\beta_n \left(n+2 \right)+n} \mathcal{G}_{\beta_n,\infty}^{(n)} (f)~.
\label{eqProposition22}
\end{equation}
If $\beta_n \neq \frac{n}{n+2}$, we can make this as small as we want, by sending $\lambda$ either to $0$ or to $\infty$. We conclude that $\mathcal{G}_{\beta_n,\infty}^{(n)}$ can only have a positive infimum, if:
\begin{equation}
\beta_n = \frac{n}{n+2}~.
\label{eqProposition23}
\end{equation}

Next, we consider example \ref{Example2}. From equations (\ref{eqPrep9}) and (\ref{eqPrep11}), we obtain:
\begin{equation}
\begin{array}{c}
\mathcal{G}_{\frac{n}{n+2},\infty}^{(n)}\left(g_c^{(n)} \right)= \left\{\frac{n}{\pi} \left(\sqrt{2}+\frac{2c}{\sqrt{c^4+1}}\right)^{n-1}\times\right. \\
\\
\left.\times \left[\frac{1}{4\sqrt{2}}\left(\frac{c^4+1}{c^2}\right) +\frac{c^3}{(c^4+1)^{3/2}}\right] \frac{c^n}{(c+1)^{2n}}\right\}^{\frac{n}{n+2}}\\
\\
\sim \left(\frac{n (\sqrt{2})^{n-6}}{\pi}\right)^{\frac{n}{n+2}} c^{\frac{n(n-2)}{n+2}} + \mathcal{O}\left( c^{\frac{(n-1)n}{n+2}}\right)~.
\end{array}
\label{eqProposition24}
\end{equation}
By sending $c\to 0^+$, we can make this arbitrarily small, whenever $n>2$.

We are left with the case $n=2$:
\begin{equation}
\mathcal{G}_{\frac{1}{2},\infty}^{(2)} (f)=\left(\frac{V(f)}{\|f\|_{\infty}^2}\right)^{\frac{1}{2}} \frac{\|f\|_{\infty}}{\|f\|_1}=\left(\frac{V(f)}{\|f\|_1^2}\right)^{\frac{1}{2}}~.
\label{eqProposition25}
\end{equation}
We consider example \ref{Example4}. From equations (\ref{eqExample43}) and (\ref{eqExample44}), we obtain:
\begin{equation}
\mathcal{G}_{\frac{1}{2},\infty}^{(2)} \left(f_{\alpha}\right)= \sqrt{\frac{5 \pi \alpha}{32}}~.
\label{eqProposition26}
\end{equation}
If $\alpha \to 0^+$, then $\mathcal{G}_{\frac{1}{2},\infty}^{(2)} \left(f_{\alpha}\right) \to 0^+$.

\subsection{Proof of Proposition \ref{Proposition2}}

Let us consider the functional 
\begin{equation}
\mathcal{G}_{\beta_{n,q},q}^{(n)} (f)= \left(\frac{V(f)}{\|f\|_{q}^2}\right)^{\beta_{n,q}} \frac{\|f\|_{q}}{\|f\|_1}~,
\label{eqProposition27}
\end{equation}
for $\beta_{n,q}>0$ and $f \in \mathcal{S}(\mathbb{R}^n)$. 

If we consider again the dilation $f_{\lambda}(x)=f(\lambda x)$, we obtain from (\ref{eqExtra1},\ref{eqExtra3}):
\begin{equation}
\mathcal{G}_{\beta_{n,q},q}^{(n)} \left(f_{\lambda} \right)=\lambda^{\beta_{n,q} \left(\frac{2n}{q}-n-2 \right)+\frac{n}{q^{\prime}} } \mathcal{G}_{\beta_{n,q},q}^{(n)} (f)~.
\label{eqProposition28}
\end{equation}
As previously, we conclude that there can be no strictly positive lower bound, unless:
\begin{equation}
\beta_{n,q} \left(\frac{2n}{q}-n-2 \right)+\frac{n}{q^{\prime}}=0 \Longleftrightarrow \beta_{n,q}=\frac{n(q-1)}{(n+2)q-2n}~.
\label{eqProposition29}
\end{equation}
The case $n=2$, $q=1$ is singular, so we treat $n=2$ separately. If $n=2$ and $q>1$, then from (\ref{eqProposition29}) $\beta_{2,q}=\frac{1}{2}$ and
\begin{equation}
\mathcal{G}_{\frac{1}{2},q}^{(2)} (f)= \left(\frac{V(f)}{\|f\|_{q}^2}\right)^{\frac{1}{2}} \frac{\|f\|_{q}}{\|f\|_1}~=\left(\frac{V(f)}{\|f\|_{1}^2}\right)^{\frac{1}{2}} ~.
\label{eqProposition30}
\end{equation}
If $n=2$ and $q=1$, we have from (\ref{eqProposition27}):
\begin{equation}
\mathcal{G}_{\beta_{2,1},1}^{(2)} (f)= \left(\frac{V(f)}{\|f\|_{1}^2}\right)^{\beta_{2,1}}  ~.
\label{eqProposition31}
\end{equation}
But from (\ref{eqProposition25},\ref{eqProposition26}), we know that we can make (\ref{eqProposition30},\ref{eqProposition31}) arbitrarily small. 

Next, we consider the case $n>2$. 

On the other hand, $\beta_{n,q}$ in (\ref{eqProposition29}) is strictly positive if and only if $q>1$ and:
\begin{equation}
(n+2)q-2n >0\Longleftrightarrow q > \frac{2n}{n+2} \geq 1~.
\label{eqProposition32}
\end{equation}
Next, we consider again example \ref{Example2}. From (\ref{eqPrep9},\ref{eqPrep11}) and the estimate (\ref{eqPrep12}), we obtain:
\begin{equation}
\begin{array}{c}
\mathcal{G}_{\frac{n(q-1)}{(n+2)q-2n},q}^{(n)} \left(g_c^{(n)} \right)=\frac{\left[V(g_c^{(n)})\right]^{\frac{n(q-1)}{(n+2)q-2n}}}{\|g_c^{(n)}\|_1 ~\|g_c^{(n)}\|_q^{\frac{(n-2)q}{(n+2)q-2n}}} \leq\\
\\
\leq \left\{\frac{n}{\pi} \left(\sqrt{2} +\frac{2c}{\sqrt{c^4+1}}\right)^{n-1} \left[\frac{1}{4\sqrt{2}}\left(\frac{c^4+1}{c^2}\right)+\frac{c^3}{(c^4+1)^{3/2}}\right]\right\}^{\frac{n(q-1)}{(n+2)q-2n}} \times \\
\\ 
\times \left(\frac{2^n q^{\frac{n}{2q}}}{c^{\frac{n}{q}-\frac{n}{2}}}\right)^{\frac{(n-2)q}{(n+2)q-2n}}\times \frac{c^{n/2}}{(c+1)^n}~.
\end{array}
\label{eqProposition32}
\end{equation}
As $c \to 0^+$, we have:
\begin{equation}
\begin{array}{c}
\mathcal{G}_{\frac{n(q-1)}{(n+2)q-2n},q}^{(n)} \left(g_c^{(n)} \right) \leq \left(\frac{n2^{(n-6)/2}}{\pi}\right)^{\frac{n(q-1)}{(n+2)q-2n}} \times \left( 2^{2q} q \right)^{\frac{n(n-2)}{2q(n+2)-4n}} c^{\frac{(q-2)n(n-2)}{(n+2)q-2n}} + \\
\\
+\text{higher order terms}
\end{array}
\label{eqProposition33}
\end{equation}
If $q>2$, this can be made arbitrarily small as $c \to 0^+$.

\subsection{Proof of Lemma \ref{Lemma1}}

We start with the $n=1$ case. Let $\theta \in (0, \pi/2)$ and denote by $\mathcal{F}_{\theta}$ the fractional Fourier transform \cite{Dias2}. For simplicity, we write $\mathcal{F}$ for the Fourier transform, instead of $\mathcal{F}_{\pi/2}$.

Recall that:
\begin{equation}
\mathcal{F}_{\alpha}\mathcal{F}_{\beta}=\mathcal{F}_{\alpha + \beta}~.
\label{eqHadSpecHad14}
\end{equation}

Since $\mathcal{F}_{\theta}$ is a linear canonical transform, we have that $\mathcal{F}_{\theta} \in \mathcal{SH}_{k(\theta)}$, for some $k(\theta)>0$, which varies with $\theta$.

Set:
\begin{equation}
A=\mathcal{F}_{\theta} \hspace{0.5 cm} \text{and} \hspace{0.5 cm} B=A \mathcal{F}= \mathcal{F}_{\theta}\mathcal{F}= \mathcal{F}_{\theta+ \pi/2}~.
\label{eqHadSpecHad15}
\end{equation}
Thus:
\begin{equation}
A \in \mathcal{SH}_{k(\theta)} \hspace{0.5 cm} \text{and} \hspace{0.5 cm} B \in \mathcal{SH}_{k(\theta+ \pi/2)}~.
\label{eqHadSpecHad16}
\end{equation}

We have:
\begin{equation}
\left\{
\begin{array}{l}
A+B=A(I +\mathcal{F})\\
(A+B)^{\ast}=(I+ \mathcal{F}^{\ast})A^{\ast}
\end{array}
\right.
\label{eqHadSpecHad17}
\end{equation}
From (\ref{eqHadSpecHad16},\ref{eqHadSpecHad17}):
\begin{equation}
\begin{array}{c}
(A+B)^{\ast}(A+B)=(I+ \mathcal{F}^{\ast})A^{\ast}A (I + \mathcal{F})=k(\theta) (I+ \mathcal{F}^{\ast})(I + \mathcal{F})=\\
\\
=k(\theta) \left(I+\mathcal{F}+ \mathcal{F}^{\ast}+\mathcal{F}^{\ast}\mathcal{F} \right)=k(\theta) \left(2I +\mathcal{F}+ \mathcal{F}^{-1}\right)~,
\end{array}
\label{eqHadSpecHad18}
\end{equation}
where we used the fact that $\mathcal{F}^{-1}=\mathcal{F}^{\ast}$.
Let $h_2(x)$ be the second Hermite function. We have:
\begin{equation}
(\mathcal{F}h_2)(x)=-h_2 (x)~, \hspace{1 cm} (\mathcal{F}^{-1} h_2)(x)=-h_2 (x)~.
\label{eqHadSpecHad19}
\end{equation}
From (\ref{eqHadSpecHad18},\ref{eqHadSpecHad19}):
\begin{equation}
\left[(A+B)^{\ast}(A+B)h_2\right](x)=k(\theta) \left[ 2 h_2 (x) -h_2 (x)-h_2 (x)\right]=0~.
\label{eqHadSpecHad20}
\end{equation}
Thus:
\begin{equation}
\frac{\|(A+B)^{\ast}(A+B)h_2\|_{\infty}}{\|h_2 \|_{\infty}}=0~,
\label{eqHadSpecHad21}
\end{equation}
which proves that $A+B \notin \mathcal{H}$. 

In dimension $n>1$ we just have to consider tensor products appropriately. Notice that, since $A$ and $B$ are special Hadamard operators, the same conclusions hold for $\mathcal{SH}$. 

\subsection{Proof of Lemma \ref{Lemma2}}

First of all notice that $I \notin \mathcal{H}$ (and $I \notin \mathcal{SH}$). Indeed, let
\begin{equation}
f(x) =\left\{
\begin{array}{l l}
\frac{1}{\sqrt{x}} & ,~ \text{if}~x \in \left. \right]0,1 \left. \right]\\
0 & ,~ \text{otherwise}
\end{array}
\right.
\label{eqHadSpecHad22}
\end{equation}
Then $\|f\|_1  < \infty$, but $\|If\|_{\infty}= \infty$. Thus $I \notin \mathcal{H}$. The same conclusion holds for $n>1$.

Note that $\mathcal{F}$, $\mathcal{F}^{-1} \in \mathcal{SH}\subset \mathcal{H}$, but
\begin{equation}
\mathcal{F} \mathcal{F}^{-1} =I  \notin \mathcal{H}~.
\label{eqHadSpecHad23}
\end{equation}

\subsection{Proof of Lemma \ref{Lemma3}}

Let $C=A$ or $B$. Then:
\begin{equation}
\|UC f\|_{\infty} \leq \|U\|_{\infty \to \infty} \|C f\|_{\infty} \lesssim  \|f\|_1~.
\label{eqHadSpecHad24}
\end{equation}
On the other hand:
\begin{equation}
(UC)^{\ast}(UC)= C^{\ast} U^{\ast}UC=\alpha C^{\ast}C~.
\label{eqHadSpecHad25}
\end{equation}
Let $C=A$. From (\ref{eqHadSpecHad25}):
\begin{equation}
\|(UC)^{\ast}(UC)f\|_{\infty}= \alpha \| A^{\ast}Af\|_{\infty} \geq \alpha k \|f\|_{\infty}~.
\label{eqHadSpecHad26}
\end{equation}
From (\ref{eqHadSpecHad24},\ref{eqHadSpecHad26}), we conclude that $UA \in \mathcal{H}_{k \alpha}$.

Alternatively, let $C=B$. From (\ref{eqHadSpecHad25}):
\begin{equation}
(UB)^{\ast}(UB)= \alpha B^{\ast}B= \alpha k I~.
\label{eqHadSpecHad27}
\end{equation}
Moreover
\begin{equation}
\|UB\|_{2 \to 2}= \sqrt{\alpha k} < \infty~,
\label{eqHadSpecHad28}
\end{equation}
and
\begin{equation}
\|(UB)^{\ast}\|_{1 \to \infty}= \|B^{\ast}U^{\ast}\|_{1 \to \infty} \leq  \|B^{\ast}\|_{1 \to \infty}  \|U^{\ast}\|_{1 \to 1}~.
\label{eqHadSpecHad29}
\end{equation}
From (\ref{eqHadSpecHad24},\ref{eqHadSpecHad27},\ref{eqHadSpecHad28},\ref{eqHadSpecHad29}) it follows that $UB \in \mathcal{SH}_{\alpha k}$.

\subsection{Example \ref{ExampleCompStepFunc}}\label{ProofExampleStepFunc}

We have:
\begin{equation}
\begin{array}{c}
\|AB\|_{\infty} \leq \|Bf\|_1= \int_{\mathbb{R}^n }\left|\sum \alpha_j (P_jf)(x) \right| dx=\\
\\
=\sum_{j \in \mathbb{N}} \alpha_j \int_{\Omega_j}|f(x)|dx \leq M \|f\|_1~.
\end{array}
\label{eqHadSpecHad30}
\end{equation}
On the other hand:
\begin{equation}
\begin{array}{c}
\left((AB)^{\ast}AB f \right)(x)=\left(B^{\ast}A^{\ast}ABf \right)(x)=\\
\\
=k (B^{\ast}B f)(x)= k \left(\sum_{j \in \mathbb{N}}\alpha_j^2 (P_jf)(x)\right)~,
\end{array}
\label{eqHadSpecHad31}
\end{equation}
where we used (\ref{eqHadSpecHad16A}).

It follows that:
\begin{equation}
\|(AB)^{\ast}AB f\|_{\infty}=k \|\sum_{j \in \mathbb{N}}\alpha_j^2 P_jf \|_{\infty} \geq km^2 \|f\|_{\infty}~.
\label{eqHadSpecHad32}
\end{equation}
From (\ref{eqHadSpecHad30}) and (\ref{eqHadSpecHad32}) we conclude that $AB \in \mathcal{H}_{km^2}$.

\subsection{Proof of Proposition \ref{PropositionInclusion1}}

Consider as in Example \ref{ExampleCompStepFunc} the operator
\begin{equation}
B= \alpha P_1+\beta P_2~,
\label{eqHadSpecHad33}
\end{equation}
where $0 < \alpha < \beta$, $\Omega \in \mathbb{R}^n$ is a measurable set, $\Omega^c$ its complement, and $(P_1f)(x)= \chi_{\Omega} f(x)$, $(P_2f)(x)= \chi_{\Omega^c} f(x)$.

If $A \in \mathcal{SH}_k$, then $AB \in \mathcal{H}_{k \alpha^2}$. On the other hand from (\ref{eqHadSpecHad31}), we have:
\begin{equation}
\left((AB)^{\ast}ABf \right)(x)= k \left(\alpha^2 (P_1f)(x)+ \beta^2 (P_2 f)(x) \right)~.
\label{eqHadSpecHad34}
\end{equation}
Since $\alpha \neq \beta$, there exists no $\gamma>0$, such that
$$
\left((AB)^{\ast}ABf \right)(x)= \gamma \left( (P_1f)(x)+ (P_2 f)(x) \right)= \gamma f(x)~.
$$
Thus, $AB \notin \mathcal{SH}$.

\subsection{Proof of Proposition \ref{PropositionInclusion2}}

Let us consider again the operator $AB$ defined in the proof of Proposition \ref{PropositionInclusion1}. We proved that $AB \in \mathcal{H}_{k \alpha^2}$ but $AB \notin \mathcal{SH}$.

On the other hand, since $A \in \mathcal{SH}_k$, we also have $A\in \mathcal{A}_{p,q}$ for all $p \in \left[1,2 \right]$ and $q \in \left[p,p^{\prime}\right]$. Thus:
\begin{equation}
\begin{array}{c}
\|ABf\|_p \|f\|_p=\|ABf\|_p \|Bf\|_p \frac{\|f\|_p}{\|Bf\|_p}\geq C_{p,q} \frac{\|ABf\|_q \|Bf\|_q \|f\|_q}{\|Bf\|_p}=\\
\\
=C_{p,q} \|ABf\|_q \|f\|_q \left(\frac{\|Bf\|_q \|f\|_p}{\|Bf\|_p \|Bf\|_q}\right)
\end{array}
\label{eqHadSpecHad35}
\end{equation}
Given $1 \leq r < \infty$, we have:
\begin{equation}
\begin{array}{c}
\|Bf\|_r= \|\alpha P_1 f+\beta P_2 f\|_r=\\
\\
=\left(\alpha^r \int_{\Omega} |f(x)|^r dx+\beta^r \int_{\Omega^c} |f(x)|^r dx\right)^{\frac{1}{r}}~.
\end{array}
\label{eqHadSpecHad36}
\end{equation}
It follows that:
\begin{equation}
\alpha \|f\|_r \leq \|Bf\|_r \leq \beta \|f\|_r~.
\label{eqHadSpecHad37}
\end{equation}
In the same fashion we can prove that
\begin{equation}
\alpha \|f\|_{\infty} \leq \|Bf\|_{\infty} \leq \beta \|f\|_{\infty}~.
\label{eqHadSpecHad38}
\end{equation}
From (\ref{eqHadSpecHad35})-(\ref{eqHadSpecHad38}):
\begin{equation}
\|ABf\|_p \|f\|_p \geq \frac{C_{p,q}\alpha}{\beta}~,
\label{eqHadSpecHad39}
\end{equation}
and $AB \in \mathcal{A}_{p,q}$ for all $p \in \left[1,2\right]$, $q \in \left[p,p^{\prime}\right]$.

\section*{Appendix}

In this appendix we present various examples, which are useful in some of the proofs.

\begin{example}\label{Example1}
Let $(h_k)_k$ denote the normalized Hermite functions
\begin{equation}
h_k(x)=\frac{(-1)^k e^{x^2 /2}}{\sqrt{2^k k! \sqrt{\pi}}} \frac{d^k}{dx^k}e^{-x^2}, ~x \in \mathbb{R},~k=0,1,2,3, \cdots
\label{eqPrep1}
\end{equation}
Hermite functions are eigenfunctions of the Fourier transform with a different normalization:
\[
\left(\mathcal{F}_{\frac{1}{2\pi}} h_k\right)(\xi)=(-i)^k h_k (\xi),
\]
where
\[
\left(\mathcal{F}_{\frac{1}{2\pi}} f\right)(\xi)= \frac{1}{\sqrt{2 \pi}} \int_{\mathbb{R}} f(x) e^{-i x \cdot \xi} dx~.
\]
Since 
\[
\left(\mathcal{F}_{\frac{1}{2\pi}} f\right)(\xi)=\frac{1}{\sqrt{2 \pi}}\widehat{f} \left(\frac{\xi}{2 \pi} \right)~,
\]
we conclude that:
\begin{equation}
\left(\mathcal{F}h_k\right)(\xi)=\widehat{h}_k(\xi)= (-i)^k \sqrt{2\pi} h_k (2 \pi \xi), ~k=0,1,2,3,\cdots
\label{eqPrep2}
\end{equation}
The $L^q$ norms of the Hermite functions have the following asymptotic behaviour as $k\to \infty$ (see Lemma 1.5.2 of \cite{Thangavelu}): 

\begin{enumerate}
\vspace{0.2 cm}
\item $\|h_k \|_q \sim k^{1/(2q)-1/4}$, if $1 \leq q <4$

\vspace{0.2 cm}
\item $\|h_k \|_q \sim k^{-1/8} \log n$, if $ q =4$

\vspace{0.2 cm}
\item $\|h_k \|_q \sim k^{-1/(6q)-1/12}$, if $4 < q \leq \infty$
\end{enumerate}
Here $a_k \sim b_k$ means $a_k = O (b_k)$ and $b_k = O (a_k)$.

In dimension $n>1$, we construct tensor products of Hermite functions for the multi-indices $\alpha =(k_1, \cdots,k_n)$, with $k_1,\cdots,k_n=0,1,2,3,\cdots$:
\begin{equation}
h_{\alpha}(x)=\left(h_{k_1}\otimes \cdots \otimes h_{k_n}\right)(x)= h_{k_1}(x_1) \cdot h_{k_2}(x_2) \cdots \cdot h_{k_n}(x_n)~,
\label{eqPrep3}
\end{equation}
for $ x=(x_1,\cdots,x_n) \in \mathbb{R}^n$.

It is also worth remarking that Hermite functions are eigenvectors of the Hamiltonian of the isotropic harmonic oscillator:
\begin{equation}
\left(-\frac{1}{2} \Delta + \frac{1}{2} |x|^2\right) h_{\alpha}(x)= E_{\alpha} h_{\alpha} (x)~,
\label{eqPrep4}
\end{equation}
where $\Delta= \sum_{j=1}^n \frac{\partial^2}{\partial x_j^2}$ is the Laplacian and the energy $E_{\alpha}$ is given by:
\begin{equation}
E_{\alpha}=|\alpha|+\frac{n}{2}~,
\label{eqPrep5}
\end{equation}
with $|\alpha|=k_1+k_2+ \cdots k_n$.

For future reference, we also remark that, in view of Plancherel's theorem and equations (\ref{eqPrep2},\ref{eqPrep4}):
\begin{equation}
\begin{array}{c}
V(h_{\alpha})=\int_{\mathbb{R}^n}|x|^2 |h_{\alpha}(x)|^2 dx=\\
\\
=\frac{1}{2} \int_{\mathbb{R}^n}|x|^2 |h_{\alpha}(x)|^2 dx +\frac{1}{2} \int_{\mathbb{R}^n}|x|^2 |h_{\alpha}(x)|^2 dx =\\
\\
=\frac{1}{2} \int_{\mathbb{R}^n}|x|^2 |h_{\alpha}(x)|^2 dx +\frac{1}{2(2\pi)^n} \int_{\mathbb{R}^n}|\xi|^2 |\widehat{h}_{\alpha}\left(\frac{\xi}{2\pi}\right)|^2 d\xi =\\
\\
=\frac{1}{2} \int_{\mathbb{R}^n}|x|^2 |h_{\alpha}(x)|^2 dx +\frac{1}{2} \int_{\mathbb{R}^n} |\nabla h_{\alpha}(x)|^2 dx=\\
\\
= \int_{\mathbb{R}^n}\overline{h_{\alpha}(x)}\left(-\frac{1}{2}\Delta h_{\alpha}(x) + \frac{1}{2}|x|^2 h_{\alpha}(x)\right)= E_{\alpha} \|h_{\alpha}\|_2^2 = E_{\alpha}~.
\end{array}
\label{eqPrep6}
\end{equation}

\end{example}

The second example consist of the tensor products of functions considered in \cite{Wigderson}.

\begin{example}\label{Example2}
Let 
\begin{equation}
g_{c}^{(1)} (x)= \frac{1}{\sqrt{c}}e^{-\frac{\pi}{c^2} x^2} + \sqrt{c}~e^{-\pi c^2 x^2}, ~
\label{eqPrep7}
\end{equation}
for $c>0$ and $x \in \mathbb{R}$. In dimension $n>1$, we consider:
\begin{equation}
g_{c}^{(n)} (x)= \left(g_{c}^{(1)}  \otimes \cdots \otimes g_{c}^{(1)} \right)(x)= \prod_{j=1}^n g_{c}^{(1)} (x_j)~,
\label{eqPrep8}
\end{equation}
where $x=(x_1, \cdots,x_n) \in \mathbb{R}^n$.

Since $g_{c}^{(1)}$ was considered in \cite{Wigderson}, we shall only briefly state the results needed. We have:
\begin{equation}
\|g_{c}^{(n)} \|_{\infty}=\|g_{c}^{(n)} \|_1=\left(\sqrt{c}+\frac{1}{\sqrt{c}}\right)^n~,
\label{eqPrep9}
\end{equation}
and  
\begin{equation}
\|g_{c}^{(n)} \|_2^2=\left(\sqrt{2}+\frac{2c}{\sqrt{c^4+1}}\right)^n~.
\label{eqPrep10}
\end{equation}

The variance of $g_{c}^{(n)}$ is:
\begin{equation}
\begin{array}{c}
V(g_{c}^{(n)} )=n \|g_{c}^{(1)} \|_2^{2(n-1)} V(g_{c}^{(1)} )=\\
\\
=\frac{n}{\pi} \left(\sqrt{2} +\frac{2c}{\sqrt{c^4+1}} \right)^{n-1}\left[\frac{1}{4\sqrt{2}}\left(\frac{c^4+1}{c^2}\right) +\frac{c^3}{\left(c^4+1\right)^{3/2}} \right]~.
\end{array}
\label{eqPrep11}
\end{equation}  
The $L^q$ norms for $2 \neq q \in \left(1,\infty\right)$ are more intricate, but for our purposes we use the following estimates obtained in \cite{Huang}:
\begin{equation}
\left(\frac{c^{1/q-1/2}+c^{1/2-1/q}}{2 q^{1/(2q)}} \right)^n < \|g_{c}^{(n)}\|_q <\left(\frac{c^{1/q-1/2}+c^{1/2-1/q}}{q^{1/(2q)}} \right)^n 
\label{eqPrep12}
\end{equation}
 
\end{example}

\begin{example}\label{Example4}
The last example is the following function in $2$ dimensions:
\begin{equation}
f_{\alpha} (x)= \left\{
\begin{array}{l l}
|x|^{\alpha-2} \sin (|x|^{\alpha}), & \text{if } (2 \pi)^{\frac{1}{\alpha}} \leq |x|\leq (3 \pi)^{\frac{1}{\alpha}}\\
&\\
0 ~,& \text{otherwise}
\end{array}
\right.,
\label{eqExample41}
\end{equation}
for $\alpha >0$ and $x=(x_1,x_2) \in \mathbb{R}^2$. Using the substitution
\begin{equation}
x_1=t^{\frac{1}{\alpha}} \cos (\theta),~x_2=t^{\frac{1}{\alpha}} \sin (\theta)~,
\label{eqExample43}
\end{equation}
for $2\pi \leq t \leq 3 \pi$, $0 \leq \theta < 2\pi$, we easily obtain:
\begin{equation}
\|f_{\alpha}\|_1= \frac{2 \pi}{\alpha} \int_{2 \pi}^{3 \pi} \sin (t) dt = \frac{4 \pi}{\alpha}~,
\label{eqExample43}
\end{equation}
and 
\begin{equation}
\begin{array}{c}
V(f_{\alpha})=\int_{\mathbb{R}^2} (x_1^2+x_2^2) |f_{\alpha}(x_1,x_2)|^2 d x_1 d x_2=\\
\\
= \frac{2\pi}{\alpha}\int_{2 \pi}^{3 \pi} t \sin^2(t) dt=\\
\\
= \frac{\pi}{\alpha}\int_{2 \pi}^{3 \pi} t \left(1- \cos(2t) \right) dt =\frac{5\pi^3}{2 \alpha}~.
\end{array}
\label{eqExample44}
\end{equation}

\end{example}

Author's addresses:

\begin{itemize}

\item Nuno Costa Dias: Escola Superior N\'autica Infante D. Henrique. Av. Eng. Bonneville Franco, 2770-058 Pa\c{c}o d'Arcos, Portugal\\
 and Grupo de F\'{\i}sica
Matem\'{a}tica, Departamento de Matem\'{a}tica, Instituto Superior T\'ecnico,
Universidade de Lisboa, Av. Rovisco Pais, 1049-001 Lisboa,
Portugal

\item Franz Luef: Department of Mathematical Sciences, NTNU Trondheim, 7041 Trondheim, Norway

\item Jo\~{a}o Nuno Prata:  Departamento de Matem\'atica, ISCTE Instituto Universit\'ario de Lisboa,
Avenida das For\c{c}as Armadas, 1649-026 Lisboa, \\
and Grupo de F\'{\i}sica
Matem\'{a}tica, Departamento de Matem\'{a}tica, Instituto Superior T\'ecnico,
Universidade de Lisboa, Av. Rovisco Pais, 1049-001 Lisboa,
Portugal

\end{itemize}


\begin{thebibliography}{999999}
\bibitem{Tang2}
Y. Tang, Heisenberg uncertainty principle and its analogues in higher dimension: via Wigderson' method. Available at: Arxiv math.FA 2312.15544.
\end{thebibliography}
\end{document}